\documentclass[10pt,letterpaper]{amsart}
\usepackage[utf8]{inputenc}

\usepackage{amsaddr}

\usepackage{amsmath}
\let\amsboldsymbol\boldsymbol
\usepackage{bm}% ändert \boldsymbol
\let\boldsymbol\amsboldsymbol
\usepackage{colortbl}
\usepackage{tabularx}
\usepackage{amssymb}
\usepackage{graphicx}
\usepackage{url}
\usepackage[colorinlistoftodos]{todonotes}
\usepackage{pgfplots} % for axis environment
\usepackage{tikz} % for tikz pictures
\tikzset{every path/.append style={line width=1pt}}

\usepackage{geometry}
\definecolor{gray}{rgb}{0.4,0.4,0.4}
\definecolor{cyan}{rgb}{0.0,0.6,0.6}
\definecolor{luh-dark-blue}{rgb}{0, 0.313, 0.608} % 00509b
\definecolor{luh-light-blue}{rgb}{0.6, 0.725, 0.847} % 99b9d8
\definecolor{darkgreen}{rgb}{0.0, 0.5, 0.0}
\definecolor{darkyellow}{rgb}{1.0, 0.88, 0.21}
\definecolor{darkred}{rgb}{0.9, 0.17, 0.31}
\definecolor{lila}{rgb}{0.57, 0.36, 0.51}
\definecolor{dark-blue}{rgb}{0.2, 0.2, 0.6}
\definecolor{Maroon}{rgb}{0.48, 0.07, 0.07}
\DeclareMathOperator{\tr}{tr}
\usepackage{lipsum}
\usepackage{amsfonts}
\usepackage{graphicx}
\usepackage{epstopdf}
\usepackage{algorithmic}
\ifpdf
  \DeclareGraphicsExtensions{.eps,.pdf,.png,.jpg}
\else
  \DeclareGraphicsExtensions{.eps}
\fi

\newtheorem*{remark}{Remark}
\newtheorem{Problem}{Problem}
\newtheorem{algorithm}{Algorithm}

\usepackage{mathtools}

\title{Robust preconditioning for a mixed formulation of phase-field fracture problems%
      }

\author{Timo Heister}
\address{Clemson University, O-110 Martin Hall, Clemson, SC\\ heister@clemson.edu}
\author{Katrin Mang}
\address{Leibniz Universit\"{a}t Hannover, Welfengarten 1, 30167 Hannover, Germany\\
  {mang@ifam.uni-hannover.de}}
\author{Thomas Wick}
\address{Leibniz Universit\"{a}t Hannover, Welfengarten 1, 30167 Hannover, Germany\\ {thomas.wick@ifam.uni-hannover.de}}

\begin{document}

\begin{abstract}
In this work, we consider fracture propagation in nearly incompressible and (fully) incompressible materials using a phase-field formulation.
We use a mixed form of the elasticity
equation to overcome volume locking effects and develop a robust, nonlinear and linear solver scheme and preconditioner for the resulting system.
The coupled variational inequality system, which is solved monolithically,  consists of three unknowns: displacements, pressure, and phase-field.
Nonlinearities due to coupling, constitutive laws, and crack irreversibility
are solved using a combined Newton algorithm for the nonlinearities in the partial differential equation and employing a primal-dual active set strategy for the crack irreverrsibility constraint.
The linear system in each Newton step is solved iteratively with a flexible generalized 
minimal residual method (GMRES).
The key contribution of this work is the development of a problem-specific preconditioner that leverages the saddle-point structure of the displacement and pressure variable.
Four numerical examples in pure solids and pressure-driven fractures are conducted on uniformly 
and locally refined meshes to investigate the robustness of the solver concerning the Poisson ratio as well as the discretization and regularization parameters. 
\end{abstract}

\maketitle

% \begin{keywords}
% phase-field fracture; 
% Schur-type preconditioning; 
% mixed finite elements;
% incompressible solids
% \end{keywords}

% \begin{AMS}
% % TW, Feb 2, 2022, 22:20 habe schon mal die AMS hinzugefuegt.
%   65F10, 65M22, 74R10, 74S05
% \end{AMS}

%%%%%%%%%%%%%%%%%%%%%%%%%%%%%%%%%%%%%%%%%%%%%%%%%%
\section{Introduction}
Phase-field fracture modeling~\cite{kuhn2010continuum,miehe2010thermodynamically} 
emerged from a variational 
formulation introduced in~\cite{FraMar98,BourFraMar00} is an attractive model approach to simulate crack propagation in solids. To date, displacement-based formulations have been used in the large majority of investigations~\cite{Bour07,miehe2010thermodynamically,kuhn2010continuum,kopanivcakova2020recursive,farrell2017linear,TANNE201880,BuOrSue10,AmGeraLoren15,BrWiBeNoRa20}.
However, considering (nearly) incompressible solids, these models are subject to locking effects, i.e., the values of the displacement field are underestimated. For this reason, mixed phase-field formulations have
been recently developed in~\cite{mang2020phase,Ma21_phd} in which 
classical ideas from (non-fractured) solids were employed by introducing 
a Lagrange multiplier for the pressure variable.
With help of a mixed form we get a stable problem formulation up to the incompressible limit~\cite{basava2020adaptive}. 
However only 
sparse direct solvers, e.g.,~\cite{saad2003iterative}, were utilized in these previous studies for solving the arising linear equation
systems.

The main purpose of the current work is to propose (for the first time)
a preconditioned iterative linear solver for solving mixed formulations of phase-field 
fracture problems. Therein, we deal with three unknowns, namely 
displacements, pressure, and phase-field, $U := (u, p, \varphi)$.
For classical, displacement-based $(u,\varphi)$ formulations iterative 
linear and multigrid methods are known. We also note that 
we consider fractures in pure solids as well as pressurized cracks 
(Sneddons's test, see e.g.,~\cite{SneddLow69,BourChuYo12,schroder2021selection}).
The reader should not confound the pressure $p$ introduced due to the mixed formulation
with the pressure $\rho$, which is imposed inside the crack region in pressurized (i.e., pressure-driven) configurations.

The first study with a clear focus on 
linear solvers is~\cite{farrell2017linear}. Therein, 
a nonlinear Gauss-Seidel scheme was proposed together with a Schur complement
based preconditioner for the linear systems. 
A parallel GMRES (generalized minimal residual) solver with diagonal preconditioner with algebraic multigrid preconditioning was developed in~\cite{heister2018parallel,heister2020pfm}.
Earlier versions were used in~\cite{heister2015primal,LeeWheWi16}, however without studying the parallel performance and scalability.
A GMRES solver with a matrix-free geometric multigrid preconditioner was 
later suggested in~\cite{jodlbauer2020matrix} with a subsequent 
parallel version in~\cite{JoLaWi20_parallel}. An overall summary 
of these developments can be found in the PhD thesis of Jodlbauer~\cite{Jo21}.
%TODO (aktuell fuer TW erledigt; aber bitte nochmal vergewissern): FURTHER LITERATURE , 1 PARAGRAPH and detailed description of existing solvers. 

A Galerkin finite element discretization yields a nonlinear system of the form $\mathcal{M}U = F$ with a $(3\times 3)$ block matrix $\mathcal{M}\in \mathbb{R}^{n\times n}, U$ and $F\in \mathbb{R}^n$. For inf-sup stability a Taylor-Hood element $Q_2^c/Q_1^c$ is used for the $(u,p)$ system. Here,
$Q_2^c$ denotes a continuous finite element space with bi-quadratic finite elements (we restrict the discussion to quadrilateral finite elements here). 
We note that computational comparisons to
stabilized low-order equal-order finite elements were undertaken in~\cite{Ma21_phd}. However,
it was found that this approach can not be recommended for the mixed phase-field fracture formulation combined with high Poisson ratios. The stabilizing terms contain $\nabla p$ with mesh-dependent coefficients leading to large gradients in the crack region.

The discretized system is nonlinear, for which we employ Newton's method as a nonlinear solver. Inside, the linear system is non-symmetric and therefore, we use a GMRES method.  The key contribution is the development of a block triangular preconditioner. Individual blocks
are approximated with inner solves using the conjugate gradient method (CG) and algebraic multigrid (AMG) from the ML package~\cite{tuminaro2000parallel,gee2006ml}. 
The mixed form of the elasticity equation has a saddle point structure, which allows to reuse spectral approximations for the inverse matrices from the Stokes problem~\cite{benzi2005numerical,boffi2008finite,DrzisgaJohnRuedeWohlmuthZulehner:2018a}.
All ingredients of the preconditioner can be parallelized and the developed code is parallel (since 
extended from pfm-cracks~\cite{heister2020pfm} with scalability tests undertaken in~\cite{heister2018parallel}). However, 
we decided to focus on the various challenges
in robustness of the $(3\times 3)$ block system. Therefore, parallel computing studies with scalability tests are outside the scope of this paper.

The main challenge in developing the preconditioner is the interaction of various model, discretization, and material
parameters to obtain a robust approach. These are 
the spatial discretization parameter $h$ and the Poisson ratio $\nu$ (related to the 
Lam\'e coefficient $\lambda$) up to the incompressible limit $\nu=0.5$, and the regularization parameter $\kappa$ and the crack bandwidth $\epsilon$. 
% TODO : Machen wir adaptive und parallel in diesem Paper: Performing computations in parallel on uniform and locally refined meshes ???
We note that the basis of this work was developed in Section 6 of the PhD thesis
of the second author~\cite{Ma21_phd} and some preliminary results were published 
in~\cite{heister2021schur}.

The outline of this paper is as follows: In Section~\ref{sec_notation}, the notation and governing equations are introduced. Next, Section~\ref{sec_num_sol}
is the main part in which we first summarize the discretization and nonlinear solver. Then, the iterative solver and a Schur-type preconditioner are derived.
Afterward, in Section~\ref{sec_results} four numerical experiments are conducted 
to substantiate the performance of our algorithmic developments.
Our work is summarized in Section~\ref{sec_conclusions}.

%%%%%%%%%%%%%%%%%%%%%%%%%%%%%%%%%%%%%%%%%%%%%%%%%%%%%%%%%%%%%%%%%%%%%%%%%%%%%%%%%%%%%%%%
%%%%%%%%%%%%%%%%%%%%%%%%%%%%%%%%%%%%%%%%%%%%%%%%%%%%%%%%%%%%%%%%%%%%%%%%%%%%%%%%%%%%%%%%

\section{Notation and governing equations}\label{notation_equations}
\label{sec_notation}
Let $\Omega$ be an open and smooth two-dimensional domain and $T:=(0,T_{\text{end}})$ is a time (i.e., loading) interval with 
the partition $0=:t_0< t_1 < \ldots < t_N := T_{\text{end}}$. The lower-dimensional crack is approximated by a phase-field indicator function $\varphi:(\Omega\times T)\to [0,1]$ with $\varphi=0$ in the crack and $\varphi=1$ in the unbroken area. The bandwidth of the zone between broken and unbroken is named $\epsilon$. Further, a displacement function is defined as $u:(\Omega\times T)\to\mathcal{R}^2$.
In the following, the scalar-valued $L^2$-product is denoted by
$(x,y) := \int_{\Omega} x \cdot y \, d\Omega,$
whereas the vector-valued $L^2$-product is described by 
$(X,Y) := \int_\Omega X : Y \, d\Omega,$
with the Frobenius product $X : Y$ of two vectors $X$ and $Y$.
We define the usual Sobolev spaces $\mathcal{V}:= H_0^1(\Omega)^2$,
$\mathcal{W}:=  H^1(\Omega)$ and a convex subset 
% TW, Feb 6, 2022: sehr, sehr gut
$\mathcal{K}:=\{\varphi\in \mathcal{W} |\ 0\leq\varphi\leq \varphi^{n-1}\leq 1\ \text{a.e. in}\ \Omega\} \subset \mathcal{W}$ and $\mathcal{U}:=L_2(\Omega)$.
Further, the degradation function is defined as $g(\varphi)=(1- \kappa)\varphi^2 + \kappa$, 
where $\kappa$ is a sufficiently small regularization parameter. The stress tensor is defined as $\sigma(u) := 2 \mu E_{\text{lin}}(u) + \lambda \text{tr} (E_{\text{lin}}(u)) I$ with a linearized strain tensor $E_{\text{lin}}(u):= \frac{1}{2}(\nabla u +\nabla u^T)$, material dependent Lam\'{e} coefficients $\lambda$ and $\mu$, and the two-dimensional identity matrix $I$. The critical energy release rate is denoted as $G_c$. Based on this notation, the pressurized phase-field fracture model in its classical form can be formulated as follows~\cite{wick2020multiphysics}:

%%%%%%%%%%%%%%%%%%%%%%%%%%%%%%%%%%%%%%%%%%%%%%%%%%%%%%%%%%%%%%%%%%%%%%%%%%%%%%%%%%%%%%%%%%%%%%%%%%%%%%%%%%%%%%%%%%%%%%%%%%%%%%%%%%%%%%%%
\begin{Problem}[Pressurized phase-field fracture]
\label{classical_form}\mbox{}\\
Let a (constant) pressure $\rho\in L^{\infty}(\Omega)$
and the initial value $\varphi(0):= \varphi^0$ be given.
Given the previous timestep data $\varphi^{n-1}:= \varphi(t_{n-1})\in \mathcal{K}$.
Find $u:= u^n \in \mathcal{V}$ and $\varphi:=
\varphi^n \in \mathcal{K}$ for loading steps
$n=1,2,\ldots, N$ with
$\{u,\varphi\} \in \mathcal{V} \times \mathcal{K}$ such that
\begin{align*}
\Bigl(g(\tilde{\varphi}) \sigma(u)\,,E_{\text{lin}}(w)\Bigr) 
+(\tilde{\varphi}^{2} \rho, \nabla \cdot w)
=&\ 0 \quad \forall w \in \mathcal{V} ,\\
 (1-\kappa) ({\varphi} \sigma(u):E_{\text{lin}}(u), \psi {-\varphi}) +  2 (\varphi  \rho \nabla \cdot  u,\psi- \varphi)\\
+  G_c  \Bigl( -\frac{1}{\epsilon} (1-\varphi,\psi-\varphi)
+ \epsilon (\nabla
\varphi, \nabla (\psi - \varphi))   \Bigr)  \geq&\  0
\quad \forall \psi \in \mathcal{K}.
\end{align*}
\end{Problem}
In the elasticity part, a linear-in-time extrapolation with $\tilde{\varphi}\coloneqq\tilde{\varphi}(\varphi^{n-1},\varphi^{n-2})$  is used in the phase-field variable $\varphi$ to obtain a convex functional~\cite{heister2015primal}. Therein,
for $\varphi^{n-2}$ at $n=1$, we set $\varphi^{-1}:= \varphi^{0}$.

Based on Problem~\ref{classical_form} and following~\cite{mang2020phase}, we introduce a pressure $p:=\lambda\; \text{tr} (E_{\text{lin}}(u))$, which is a Lagrange multiplier. As mentioned in the introduction, 
$p$ and the crack pressure $\rho$ (see for instance~\cite{wick2020multiphysics} and therein
itself denoted as $p$) should not be mixed up.

\begin{Problem}[Pressurized phase-field fracture in mixed form]
\label{final_form}\mbox{}\\
Let $\rho\in L^{\infty}(\Omega)$ be given and the initial value $\varphi(0):= \varphi^0$ be given.
Given the previous time step data $\varphi^{n-1} \in \mathcal{K}$.
Find $u:= u^n \in \mathcal{V}$, $p:=p^n \in \mathcal{U}$ and $\varphi:=
\varphi^n \in \mathcal{K}$ for loading steps
$n=1,2,\ldots, N$ with
$U:=\{u,p,\varphi\} \in \mathcal{V} \times \mathcal{U} \times \mathcal{K}$ such that
\begin{align*}%\label{eq_mixed}
\begin{aligned}
\Bigl(g(\tilde{\varphi}) \sigma(u,p)\,,E_{\text{lin}}(w)\Bigr) 
+({\tilde{\varphi}}^{2} \rho, \nabla \cdot w)
=&\ 0\quad \forall w \in \mathcal{V},\\
g(\tilde{\varphi}) (\nabla \cdot u, q) - (\frac{1}{\lambda} p,q)
=&\ 0 \quad \forall q\in \mathcal{U},\\
 (1-\kappa) ({\varphi} \sigma(u,p):E_{\text{lin}}(u), \psi {-\varphi}) 
 +  2 (\varphi  \rho \nabla \cdot  u,\psi- \varphi)\\
 +  G_c  \Bigl( -\frac{1}{\epsilon} (1-\varphi,\psi-\varphi) + \epsilon (\nabla \varphi, \nabla (\psi - \varphi))   \Bigr)  \geq&\  0
\quad \forall \psi \in \mathcal{K},
\end{aligned}
\end{align*}
% \twick{TODO: Wann tilde varphi und wann nur varphi bei $\rho$? Dementsprechend muss dann 
% Problem 3 und auch die Newton-Matrix angepasst werden. Die Newton-Matrix ist aktuell
% in den $\rho$-Termen falsch.} stimmt nun
where the stress tensor is defined as $\sigma(u,p) := 2 \mu E_{\text{lin}}(u) + p I$.\\
\end{Problem}
\begin{remark}
It is clear that by setting $\rho = 0$, we obtain a phase-field formulation for 
fracture in pure solids. With this, we can investigate our preconditioner for 
both situations, namely fracture in solids and pressurized cracks.
\end{remark}

 %%%%%%%%%%%%%%%%%%%%%%%%%%%%%%%%%%%%%%%%%%%%%%%%%%%%%%%%%%%%%%%%%%%%%%%%%
 %%%%%%%%%%%%%%%%%%%%%%%%%%%%%%%%%%%%%%%%%%%%%%%%%%%%%%%%%%%%%%%%%%%%%%%%%
 \section{Discretization and numerical solution}
\label{sec_num_sol}
For the spatial discretization of Problem \ref{final_form}, we employ a Galerkin finite element method in each incremental step, where the domain $\Omega$ is partitioned into quadrilaterals
\cite{Cia87} with the discrete spaces $\mathcal{V}_h, \mathcal{U}_h$, and the convex set $\mathcal{K}_h\subset \mathcal{W}_h$.
To fulfill a discrete inf-sup condition, stable Taylor-Hood elements with 
continuous bi-quadratic shape functions ($Q_2^c$) for the displacement field $u$ and bilinear shape functions ($Q_1^c$)
for the pressure variable $p$ and the phase-field variable $\varphi$ are used as in~\cite{mang2020phase}.

%%%%%%%%%%%%%%%%%%%%%%%%%%%%%%%%%%%%%%%%%%%%%%%%%%%%%%%%%%%%%%%%%%%%%%%%%%%%%%%%%%%%%%%%%

\subsection{Nonlinear solver}
The nonlinear solution algorithm is based on a combined method. First,
nonlinearities arising from the PDE (partial differential equation) are treated with a standard line-search assisted Newton scheme. 
The crack irreversibility is handled with a primal-dual active set method. The combination of both techniques 
yields one single nonlinear Newton iteration; see~\cite{heister2015primal} for further details.

% \twick{INFO: in diesem Abschnitt gingen einige Indizes drunter und drueber und 
% auch die Energieform $E_{\epsilon}$ war gar nicht definitiert. Ich wuerde letzteres 
% einfach weglassen. Die Indizes habe ich hoffentlich alle korrigiert. Aber bitte 
% insbesondere nochmal alle $h$, $n$ und $k$ pruefen.}
%%%%%%%%%%%%%%%%%%%%%%%%%%%%%%%%%%%%%%%%%%%%%%%%%%%%%

\begin{Problem}[Discretized pressurized phase-field fracture in mixed form]\label{DiscrMixed}\mbox{}\\
Define $u_h^n:=u_h(t_n), p_h^n:=p_h(t_n)$ 
and $\varphi_h^n:= \varphi_h(t_n)$ at the loading step $t_n$.
Let $\tilde{\varphi}_h := \tilde{\varphi}_h (\varphi_h^{n-1}, \varphi_h^{n-2})$ be the discrete
linear-in-time extrapolation.
Find $U_h^n:=(u_h^n,p_h^n,\varphi_h^n)\in \mathcal{V}_h\times \mathcal{U}_h \times \mathcal{K}_h$ for all loading steps $n=1,2,..,N$ such that 
\begin{align*}
     A(U_h^n)(\Psi_h - \Phi_h^n)=&\ A_1(U_h^n)(w_h)+A_2(U_h^n)(q_h)+A_3(U_h^n)(\psi_h - \varphi_h^n) \geq 0
     \end{align*}
     with $\Phi_h^n = (0,0,\varphi_h^n)$ and
for all $\Psi_h:= (w_h,q_h,\psi_h)\in \mathcal{V}_h\times \mathcal{U}_h \times \mathcal{K}_h$, and 
where
     \begin{align*}
     \ A_1(U_h^n)(w_h)=&\ g(\tilde{\varphi}_h) \left(\sigma(u_h^n,p_h^n), E_{\text{lin}}(w_h) \right) +(\tilde{\varphi}_h^{2} \rho, \nabla \cdot w_h),\\
     A_2(U_h^n)(q_h) =&\ g(\tilde{\varphi}_h)\left(\nabla \cdot u_h^n,q_h\right) - \frac{1}{\lambda} (p_h^n,q_h),\\
     A_3(U_h^n)(\psi_h - \varphi_h^n)=&\ (1- \kappa)\left(\varphi_h^n \sigma(u_h^n,p_h^n) : E_{\text{lin}}(u_h^n),\psi_h-\varphi_h^n\right)\\
     +&\  2 (\varphi_h^n \rho \nabla \cdot  u_h^n,\psi_h - \varphi_h^n)\\
  +&\ G_c\left( (-\frac{1}{\epsilon} (1-\varphi_h^n),\psi_h-\varphi_h^n)
  + \epsilon \left(\nabla \varphi_h^n,\nabla(\psi_h-\varphi_h^n) \right)\right),
\end{align*}
where $\sigma(u_h^n,p_h^n) := 2 \mu E_{\text{lin}}(u_h^n) + p_h^n I$.
\end{Problem}
In order to treat the inequality constraint in $\mathcal{K}_h$, we employ a primal-dual 
active set method as explained in~\cite{heister2015primal} and use the 
function space $\mathcal{W}_h$ for approximating $\varphi$.
Then, at each loading step $n$, we have the following Newton iteration indexed by $k$. 
We set as initial guess $U_h^{n,0}:= U_h^{n-1}$ and iterate for $k=1,2,3,\ldots$:
 \begin{align*}
     \nabla A(U_h^{k})(\delta U_h^{n,k},\Psi)=-A(U_h^{n,k})(\Psi)  \quad \forall \Psi\in  \mathcal{V}_h\times \mathcal{U}_h\times \mathcal{W}_h.
 \end{align*}
 
 The directional derivative $\nabla A(U_h^{n,k})(\delta U_h^k,\Psi)$ in direction $\delta U_h^{k}$ 
 for $U_h^{n,k} \in \mathcal{V}_h\times \mathcal{U}_h \times \mathcal{W}_h$ is given by
 \begin{align*}
 \begin{aligned}
     \nabla A(U_h^{n,k})(\delta U_h^k&,\Psi) = g(\tilde{\varphi}_h^n) \left(\sigma(\delta u_h^k,\delta p_h^k),E_{\text{lin}}(w_h)\right)\\ 
     +&\ g(\tilde{\varphi}_h^n) \left(\nabla \cdot \delta u_h^k,q_h\right) 
     - \frac{1}{\lambda} (\delta p_h^k,q_h)\\
     +&\  (1- \kappa)\left(\varphi_h^{n,k} 2\mu 
     (E_{\text{lin}}(\delta u_h^k): E_{\text{lin}}(u_h^{n,k})+ E_{\text{lin}}(u_h^{n,k}):E_{\text{lin}}(\delta u_h^k)),\psi_h\right)\\
     +&\  2 (\varphi_h^{n,k} \rho\nabla \cdot \delta u_h^k,\psi_h)
     + (1-\kappa)(\varphi_h^{n,k}\delta p_h^k I : E_{\text{lin}}(u_h^{n,k}),\psi_h)\\
     +&\ (1- \kappa)\left(\delta\varphi_h^k \sigma(u_h^{n,k},p_h^{n,k}):E_{\text{lin}}(u_h^{n,k}),\psi_h\right)
     + 2(\delta \varphi_h^{k}\rho\nabla \cdot u_h^{n,k},\psi_h)\\
  +&\  G_c \left(\frac{1}{\epsilon} \delta \varphi_h^{k},\psi_h\right)
  + G_c \epsilon \left(\nabla \delta \varphi_h^{k},\nabla\psi_h \right).
  \end{aligned}
 \end{align*}

%%%%%%%%%%%%%%%%%%%%%%%%%%%%%%%%%%%%%%%%%%%%%%%%%%%%%%%%%%%%%%%%%%%%%%%%%%%%%%%%%%%%%%%%%

\subsection{Linear solution and Schur-type preconditioning}\label{SchurPrecond}
For the arising linear systems $\mathcal{M}\delta U=F$ inside Newton's method, a GMRES method is used,
which is right-preconditioned~\cite{saad2003iterative} with a Schur-type preconditioner $P^{-1}$. 
As usual, the goal when developing $P^{-1}$
is to have the eigenvalues of $(\mathcal{M} P^{-1})$ be independent of discretization, regularization parameters and coefficients of the problem.

%%%%%%%%%%%%%%%%%%%%%%%%%%%%%%%%%%%%%%%%%%%%%%%%%%%%%%%%%%%%%%%%%%%%%%%%%%%%%%%%%%%%%%%%%
\subsubsection{Preconditioning the $(3\times 3)$ linear system}\label{linear_solver_sec}

The system matrix $\mathcal{M}_{\text{mixed}}$ of the mixed phase-field fracture from the modified mixed problem formulation has the following block structure~\cite{heister2021schur}:
\begin{align}
\mathcal{M}_{\text{mixed}}=
  \begin{pmatrix}
     M^{uu} & M^{u p}  & M^{u \varphi}\\
   M^{p u} & M^{p p} & M^{p \varphi}\\
    M^{\varphi u} & M^{\varphi p} & M^{\varphi \varphi}
  \end{pmatrix}
  =
  \begin{pmatrix}
    g(\tilde{\varphi}) A_u & g(\tilde{\varphi}) B^T & 0\\
    g(\tilde{\varphi}) B & -\frac{1}{\lambda} M_p & 0\\
    E & F & L 
  \end{pmatrix},\label{M}
  \end{align}
  %\twick{TODO: Ich finde das zu kurz: wir sollten dem Leser helfen und formal mathematisch hinschreiben
  %welcher Block welchem Term in der Jacobi-Matrix entspricht. Oder halt mit Seitennummer genau auf 
  %Deine Dissertation verweisen.} Reicht Dir der Verweis auf die Seite in der Diss?
  where block $A_u$ is the mass matrix of the displacements, $B$ and $B^T$ are symmetric off-diagonal blocks coupling $u$ and $p$, and $M_p$ is the mass matrix of the pressure variable. The blocks $E$, $F$ and $L$ from Equation~\eqref{M} consist of the entries from the phase-field equation, where $L$ is Laplacian-like. For the entry-wise definition of the blocks, we refer to~\cite[page 169]{Ma21_phd}.
  
  A typical block factorization of the system matrix yields the preconditioner 
  (details can be found in~\cite[Chapter 6]{Ma21_phd})
 \begin{align*}
  P_{\text{mixed}}=
%   \begin{pmatrix}
%      P^{uu} & P^{u p}  & P^{u \varphi}\\
%   P^{p u} & P^{p p} & P^{p \varphi}\\
%     P^{\varphi u} & P^{\varphi p} & P^{\varphi \varphi}
%   \end{pmatrix}
%  =
   \begin{pmatrix}
    g(\tilde{\varphi}) A_u & g(\tilde{\varphi}) B^T & 0\\
    0 & S & 0\\
    0 & 0 & L 
  \end{pmatrix},
  \end{align*}
  for $\mathcal{M}_{\text{mixed}}$,
  where $S$ is the Schur complement block defined as
    \begin{align*}
      S = -\frac{1}{\lambda}  M_p - g(\tilde{\varphi}) B^T \cdot [g(\tilde{\varphi}) A_u]^{-1} \cdot g(\tilde{\varphi}) B.
  \end{align*}
  It is not feasible to construct $S^{-1}$ or even $S$ exactly, as this would
  result in a dense matrix. This means that exact evaluation of $P_{\text{mixed}}^{-1}$ is
  also not a feasible option, but it helps us to design
  an appropriate preconditioner by approximating the action of $P_{\text{mixed}}^{-1}$
  with an operator $\hat{P}^{-1}_{\text{mixed}}$ defined below.
  Note that all eigenvalues of $\mathcal{M}_{\text{mixed}} P_{\text{mixed}}^{-1}$ 
  are equal to one and GMRES would converge in at most two iterations~\cite{Murphy2000,benzi2005numerical}.
  
  Without considering the last row and column of $\mathcal{M}_{\text{mixed}}$ and 
  $P_{\text{mixed}}^{-1}$ (the phase-field), this is a typical saddle-point problem with a penalty term, where block triangular preconditioners are a common choice~\cite{benzi2005numerical}, first considered by Bramble and Pasciak in 1988~\cite{bramble1988preconditioning}, and frequently used for Stokes-type problems~\cite{clevenger2021comparison} and the Oseen equations~\cite{klawonn1999block}, where mesh-independent convergence can be observed. 
  
  To be able to efficiently apply $P_{\text{mixed}}^{-1}$, we require approximations of
  the inverses of the Laplacian-like matrices $L$, of $g(\tilde{\varphi}) A_u$, and of the Schur complement matrix $S$. With spectrally equivalent approximations, this would
  result in an optimal preconditioner~\cite{benzi2005numerical} yielding an eigenvalue distribution independent of mesh size $h$ and other problem parameters and therefore constant GMRES iterations numbers independent of mesh size and problem parameters.
  Since multigrid methods allow for mesh-independent convergence~\cite{hackbusch2013multi}, algebraic or geometric multigrid methods are the method of choice.

%  \twick{TODO: wird `remains unclear' spaeter in %diesem Paper aufgeloest oder ist 
%  das insgesamt eine offene Frage?} nicht unklar. Mehrgitter erlauben erstmal Gitterunabhaengigkeit. Auf die Abhaengigkeit von AMG von kappa gehen wir später ein. Würde ich hier nicht anführen.
  
  The approximation of the inverse of $S$ turns out to be more
  challenging. It is well-known, e.g.,~\cite{verfurth1984error}, for inf-sup stable discretizations of the linear elasticity problem, the Schur-complement is spectrally equivalent to the mass matrix.
  In our case, for $\lambda\to\infty$, $g(\tilde{\varphi})$ acts like a varying viscosity. It is common to scale
  the mass matrix with the inverse of the viscosity for Stokes interface problems ~\cite{OlshanskiiReuskenStokesInterface} or variable viscosity Stokes problems, e.g.,~\cite{Grinevich2009,May2015HetStokes},
  which yields
  \begin{align*}
      \hat{S}^{-1}\coloneqq -\left(\left(\frac{1}{\lambda}  + \frac{g(\tilde{\varphi})}{2\mu}\right) M_p\right)^{-1},
  \end{align*}
  as an approximation of the inverse of $S$ in our situation. 
  Under sufficient regularity, $\kappa>0$ and if the coefficient can be assumed to be constant,  $\hat{S}^{-1}$ is 
  spectrally equivalent to $S^{-1}$~\cite{OlshanskiiReuskenStokesInterface}.
  For the incompressible limit $\nu=0.5$, the Schur complement approximation becomes
 \begin{align*}
     \hat{S}^{-1} = -\left(\frac{g(\tilde{\varphi})}{2\mu}M_p\right)^{-1}.
 \end{align*}
  
  \begin{remark}[Differences to Stokes-type problems]\label{remark_diff}
  Commonly, this Schur complement approach is used for Stokes-type problems and incompressible fluid dynamics, see, e.g.~\cite{elman2014finite}.
  Even if the elasticity part of the considered phase-field fracture problem has a similar saddle-point structure, aside from the phase-field function, material and regularization parameters complicate the situation: $\lambda\to\infty$ leads to a purely $\kappa$-dependent block $\hat{S}^{-1}$, and $\kappa\to 0$ increases the condition number of the block $(g(\tilde{\varphi}) A_u)^{-1}$ in the crack, where $\varphi=0$.   
  While the approximation of $\hat{S}^{-1}$ is spectrally equivalent with respect to the mesh size, it is not robust with respect to large viscosity variations, or in our case minimum and maximum
  value of $g(\tilde{\varphi})$ throughout the domain. 
  For the Stokes interface problem with a viscosity jump with single interface, the scaled mass matrix is spectrally equivalent independent of the magnitude of the jump~\cite{OlshanskiiReuskenStokesInterface}, which is the case in our situation.
  This will be visible in
   Section~\ref{sec_results}. 
   We hypothesize that a better Schur complement could be a weighted BFBT preconditioner presented in~\cite{rudi2017weighted}, but a thorough investigation is future work.
  \end{remark}

%%%%%%%%%%%%%%%%%%%%%%%%%%%%%%%%%%%%%%%%%%%%
\subsubsection{Preconditioning algorithm}\label{subsec_algo_precond}

As discussed above, the evaluation of the preconditioner
\begin{align*}
\hat{P}^{-1}_{\text{mixed}}=&\
  \begin{pmatrix}
    (g(\tilde{\varphi})A_u)^{-1} & -A_u^{-1} B^T \hat{S}^{-1} & 0\\
    0 & \hat{S}^{-1} & 0\\
    0 & 0 & L^{-1}
  \end{pmatrix}
\end{align*}
requires efficient approximations to the exact inverses of $A_u$, $\hat{S}$, and $L$.
Iterative solvers like GMRES of course only require the result of a matrix-vector
product with the preconditioner $\hat{P}^{-1}_{\text{mixed}}$, 
see~\cite{saad2003iterative} and inside our basis software deal.II~\cite{arndt2021deal}, see~\cite{kronbichler2012high}.

First, we approximate $L^{-1}$ by a single $V$-cycle of algebraic multigrid (AMG).
Second, for $(g(\tilde{\varphi})A_u)^{-1}$ we use an inner Conjugate Gradient (CG) solve, which, in turn, is preconditioned by one $V$-cycle of algebraic multigrid.
Finally, the action of $\hat{S}^{-1}$ is either done using a single $V$-cycle of AMG or, in Figures~\ref{Sneddon2dKappSmallCgforS} and~\ref{Sneddon2dlayeredAdaptiveCgforS}, using CG preconditioned by AMG.

With this, the matrix-vector product $\hat{P}^{-1}_{\text{mixed}}\vec{x}$ with $\vec{x}=(x_u,x_p,x_{\varphi})^T$ given as
\begin{align*}
%\begin{aligned}
\hat{P}^{-1}_{\text{mixed}}\vec{x}=&\
%   \begin{pmatrix}
%     (g(\tilde{\varphi})A_u)^{-1} & -(g(\tilde{\varphi})A_u)^{-1} g(\tilde{\varphi}) B^T \hat{S}^{-1} & 0\\
%     0 & \hat{S}^{-1} & 0\\
%     0 & 0 & L^{-1}
%   \end{pmatrix}
%   \begin{pmatrix}
%     x_u\\
%     x_p\\
%     x_{\varphi}
%   \end{pmatrix}\\
%   =&\
  \begin{pmatrix}
   (g(\tilde{\varphi})A_u)^{-1} x_u -(g(\tilde{\varphi})A_u)^{-1} g(\tilde{\varphi}) B^T \hat{S}^{-1} x_p\\
    \hat{S}^{-1} x_p\\
    L^{-1} x_{\varphi}
  \end{pmatrix},
 % \end{aligned}
  \end{align*}
is built up step by step. 
In deal.II~\cite{dealII92,arndt2021deal}, the preconditioners given to solver classes need a \texttt{vmult()} member function~\cite{step20}. Then, our final algorithm is designed as follows:
%.The matrix-vector product to build $\hat{P}^{-1}_{\text{mixed}}\vec{x}$, is done in the member function \texttt{vmult()} in a solver class.
%https://www.dealii.org/current/doxygen/deal.II/step_20.html
%%%%%%%%%%%%%%%%%%%%%%%%%%%%%%%%%%%%%%%%%%%%%%%%%%%%%%%
\begin{algorithm}{Evaluation of $\hat{P}^{-1}_{\text{mixed}}\vec{x}$:}\label{alg_p}\mbox{}
\begin{enumerate}
    \item Approximate $\hat{S}^{-1}$ via AMG and compute $q\coloneqq\hat{S}^{-1} x_p$;
    \item Compute $r\coloneqq x_u-g(\tilde{\varphi})B^T q$;
    \item Approximate $(g(\tilde{\varphi})A_u)^{-1}$ via CG preconditioned with AMG and compute $s\coloneqq(g(\tilde{\varphi})A_u)^{-1} r$;
    \item Approximate $L^{-1}$ via AMG and compute $t\coloneqq L^{-1}x_{\varphi}$;
    \item Return the result $(s,q,t)^T$.
\end{enumerate}
\end{algorithm}

%%%%%%%%%%%%%%%%%%%%%%%%%%%%%%%%%%%%%%%%%%%%%%%%%%%%%%%%%%%%%%%%%%%%%%%%%%%%%%%%%%%%%%%%%
%%%%%%%%%%%%%%%%%%%%%%%%%%%%%%%%%%%%%%%%%%%%%%%%%%%%%%%%%%%%%%%%%%%%%%%%%%%%%%%%%%%%%%%%%
\section{Numerical tests}
 \label{sec_results}
In this section, we consider four different numerical experiments to substantiate 
our algorithmic developments and to investigate the performance of the 
nonlinear solver, linear solver and preconditioner.

\subsection{Test cases and presentation of our results}
To facilitate the readability of the tables from the next sections, we give an overview, how to read them. For the four tests, we conduct numerical studies with different emphases: we investigate robustness in $h$, $\kappa$, $\lambda$, $\epsilon$, we use different models (`primal' from Problem~\ref{classical_form} versus `mixed' from Problem~\ref{final_form}) and different finite element discretizations.
 In the top row of each table, we summarize the key aspect of the current numerical study: the name of the example, the observed task, the modeling, and -- if required -- further test-specific settings. The white rows in the tables correspond to results based on the primal 
 phase-field fracture model (solved with pfm-cracks~\cite{heister2020pfm}) or to reference values. The colored rows belong to computations based on the mixed model and $Q_2^c Q_1^c Q_1^c$ for $\nu=0.2$ (yellow), $\nu=0.4999$ (blue) and $\nu=0.5$ (red). 
 A more saturated shading denotes a finer mesh size.
 
 The four test configurations with attributes are given in the following:
 \begin{itemize}
     \item Section~\ref{sec_hbs}: a hanging block with an initial slit for $\nu=0.2,0.4999$ and $0.5$, uniform mesh refinement, mixed ($Q_2^c Q_1^c Q_1^c$) versus primal ($Q_2^c Q_1^c$), $\epsilon$ fixed and $\epsilon=2\,h$, $\kappa=10^{-2}$;
     \item Section~\ref{sec_sneddon}: Sneddon's test~\cite{sneddon1946distribution,SneddLow69} for $\nu=0.2,0.4999$ and $0.5$, uniform mesh refinement, mixed ($Q_2^c Q_1^c Q_1^c$), $\epsilon$ fixed and $\epsilon=2\,h$, $\kappa=10^{-2}, 10^{-8}$;
     \item Section~\ref{sec_sneddon_layered}: Sneddon's test layered~\cite{basava2020adaptive} for $\nu=0.2,0.4999$ and $0.5$ in the inner domain, adaptive mesh refinement (geometric), mixed ($Q_2^c Q_1^c Q_1^c$), $\epsilon=h$, $\kappa=10^{-2}, 10^{-8}$;
     \item Section~\ref{sec_sent}: single-edge notched tension test for $\nu=0.3,0.45$, and $0.49$, adaptive mesh refinement (predictor-corrector scheme), mixed ($Q_2^c Q_1^c Q_1^c$), $\epsilon=4\,h$, $\kappa=10^{-8}$.
    % \item Sneddon's test layered in 3d for $\nu=0.5$ in the inner domain (uniform and adaptive meshes)
 \end{itemize}
 With the help of numerical studies, we investigate the robustness of the new Schur-type preconditioner via evaluating the required number of linear iterations for different mesh sizes, Poisson ratios, $\kappa$, and different finite element discretizations. Besides, we discuss challenges and point out difficulties.
 
 %%%%%%%%%%%%%%%%%%%%%%%%%%%%%%%%%%%%%%%%%%%%%%%%%%%%%%%%%%%%%%%%%%%%%%%%%%%%%%%%%%%%
\subsection{Implementation details}

The software developed for this paper is a major extension built upon pfm-cracks~\cite{heister2018parallel,heister2020pfm}, 
which is an open-source code available
at \url{https://github.com/tjhei/cracks}. This project is built on the finite element library deal.II~\cite{dealII92}, which offers scalable parallel algorithms for finite element computations.
The deal.II library in turn uses functionality from other libraries such as Trilinos~\cite{heroux2005overview,trilinos} for linear algebra, including the Trilinos ML AMG preconditioner~\cite{tuminaro2000parallel,gee2006ml}.
The GMRES stopping criterion is a relative tolerance of $10^{-5}$. CG uses a relative tolerance of $10^{-6}$ for the inner solves with a maximum of 200 iterations.
The Newton iteration stops when an absolute tolerance of $10^{-7}$ is reached.
We use four CPUs on a single machine with four Intel E7 v3 CPUs for all computations.

%%%%%%%%%%%%%%%%%%%%%%%%%%%%%%%%%%%%%%%%%%%%%%%%%%%%%%%%%%%%%%%%%%%%%%%%%%%%%%%%%%%%

 \subsection{Hanging block with initial slit}\label{sec_hbs}

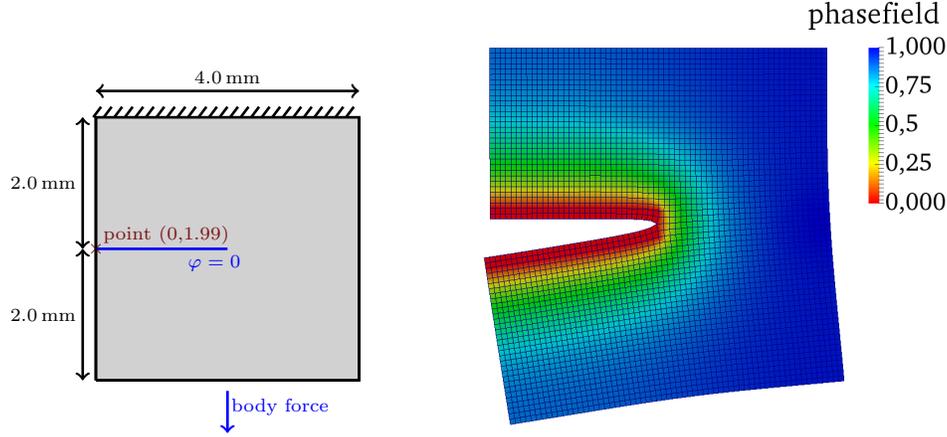
\begin{figure}[htbp!]
 \centering
 \scriptsize
 \begin{tikzpicture}[xscale=0.7,yscale=0.7]
\draw[fill=gray!30] (0,0) -- (0,5) -- (5,5) -- (5,0) -- (0,0);
\draw (0,0) -- (5,0);
\draw (5,0) -- (5,5);
\draw[blue] (0,2.5) -- (2.5,2.5);
\node[Maroon] at (0,2.5) {$\times$};
\node[Maroon] at (1.34,2.76) {point (0,1.99)}; 
%\draw[Maroon] (2.45,2.5) -- (2.5,2.5);
\node[blue] at (2.25,2.22) {$\varphi=0$};
%\draw[lu, thick] (5,5) -- (0,5);
\draw[-] (0,5) -- (0,0);
\draw[<->] (-0.25,0) -- (-0.25,2.5);
\draw[<->] (-0.25,2.5) -- (-0.25,5);
%\node at (2.5,5.2) {$\Gamma_{\text{top}}$};
%\node at (0.44,5.25) {$u_x$};
\node[blue] at (3.5,-0.5) {body force};
\draw[->,blue] (2.5,-0.2) -- (2.5,-1.0);
\node at (-1,1.25) {$2.0\,\mathrm{mm}$};
\node at (2.5,5.75) {$4.0\,\mathrm{mm}$};
\node at (-1,3.75) {$2.0\,\mathrm{mm}$};
\draw[<->] (0,5.5) -- (5,5.5);
%\node at (-0.65,2.5) {$4$};
\draw (0.1,5.2) -- (-0.05,5);
\draw (0.3,5.2) -- (0.15,5);
\draw (0.5,5.2) -- (0.35,5);
\draw (0.7,5.2) -- (0.55,5);
\draw (0.9,5.2) -- (0.75,5);
\draw (1.1,5.2) -- (0.95,5);
\draw (1.3,5.2) -- (1.15,5);
\draw (1.5,5.2) -- (1.35,5);
\draw (1.7,5.2) -- (1.55,5);
\draw (1.9,5.2) -- (1.75,5);
\draw (2.1,5.2) -- (1.95,5);
\draw (2.3,5.2) -- (2.15,5);
\draw (2.5,5.2) -- (2.35,5);
\draw (2.7,5.2) -- (2.55,5);
\draw (2.9,5.2) -- (2.75,5);
\draw (3.1,5.2) -- (2.95,5);
\draw (3.3,5.2) -- (3.15,5);
\draw (3.5,5.2) -- (3.35,5);
\draw (3.7,5.2) -- (3.55,5);
\draw (3.9,5.2) -- (3.75,5);
\draw (4.1,5.2) -- (3.95,5);
\draw (4.3,5.2) -- (4.15,5);
\draw (4.5,5.2) -- (4.35,5);
\draw (4.7,5.2) -- (4.55,5);
\draw (4.9,5.2) -- (4.75,5);
 \end{tikzpicture}
 \hspace{1.0cm}
 \includegraphics[width=0.4\textwidth]{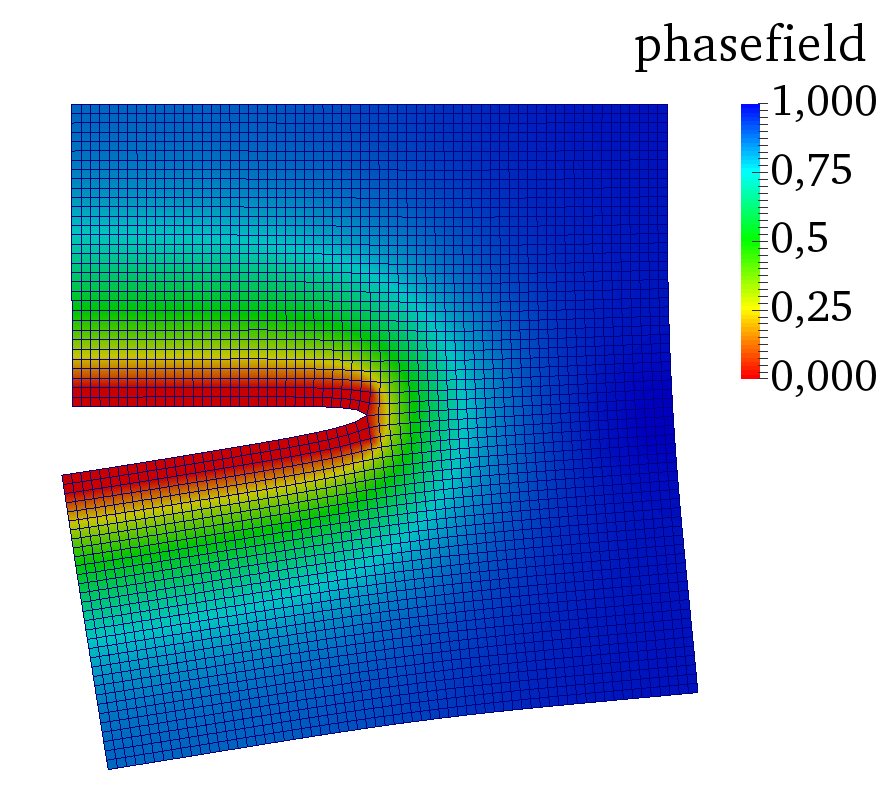}
\caption{Left: geometry and boundary conditions of a hanging block with a slit. Right: deformed geometry with phase-field solution; $41\,924$ degrees of freedom (DoFs). In the geometrically pre-refined slit we interpolate $\varphi=0$ as an initial condition. No propagating crack. Applied force $f=-8.0\cdot10^{-7}\,\mathrm{N/mm^2}$.}\label{hbs_geo}
\end{figure}

 As a first test configuration, we consider a hanging block test with an initial geometrical slit of length $2.0\,\mathrm{mm}$ with an interpolated initial condition $\varphi=0$ in the crack; see Figure~\ref{hbs_geo}. The force acting on the hanging block is reduced to $f=-8.0\cdot10^{-7}\,\mathrm{N/mm^2}$. 
In Figure~\ref{hbs_geo} on the right, the solution of the phase-field function is given on the deformed block for $\nu=0.2$ on a uniform refined mesh with $41\,924$ degrees of freedom (DoFs). We evaluate the displacement in the $y$-direction in a certain point $(0,1.99)$ on the lower opening crack lip.

\begin{table}[htbp!]
\tiny
\centering
\caption{A hanging block with an initial slit for $\nu=0.2$ and $\nu=0.4999$, uniform mesh refinement, mixed ($Q_2^c Q_1^c Q_1^c$) versus primal ($Q_2^c Q_1^c$), $\epsilon$ fixed, $\kappa=10^{-2}$. Problem size in $\#$DoFs, average number of GMRES iterations ($\varnothing$lin) per Newton/active set (AS) step, the average number of inner CG iterations ($\varnothing$CG) per linear iteration, number of Newton/AS ($\#$AS), and goal functional displacement in a certain point ($u_y(0,1.99)$). Applied force $f=-8\cdot 10^{-7}\,\mathrm{N/mm^2}$.
}\label{HangingBlockSlitTable}
\begin{tabular}{|l|l|c|c|c|r|r|r|r|l|} \hline
 \multicolumn{10}{|c|}{Hanging block slit: robustness in $h$ and $\lambda$; mixed versus primal; $\kappa=10^{-2}$}\\ \hline\hline
 \rowcolor{gray!50} \multicolumn{1}{|c|}{model} & \multicolumn{1}{|c|}{FE} & \multicolumn{1}{|c|}{$\nu$} & \multicolumn{1}{|c|}{$h$ }& \multicolumn{1}{|c|}{$\epsilon$} & \multicolumn{1}{|c|}{$\#$DoFs} & \multicolumn{1}{|c|}{$\varnothing$lin} & \multicolumn{1}{|c|}{$\varnothing$CG} & \multicolumn{1}{|c|}{$\#$AS} & \multicolumn{1}{|c|}{$u_y(0,1.99)$}  \\
 \hline
  \hline
 \rowcolor{yellow!10}mixed & $Q_2^c Q_1^c Q_1^c$ & 0.2 & $0.353$ & $0.707$ & $2\,804$ & 4 & 24 & 3 &  -0.3871  \\ \hline %3
  \rowcolor{yellow!20}mixed & $Q_2^c Q_1^c Q_1^c$ & 0.2 & $0.176$  & $0.707$ & $10\,724$  &4  & 25 & 3 & -0.5189\\ \hline %4
    \rowcolor{yellow!30}mixed & $Q_2^c Q_1^c Q_1^c$& 0.2 & $0.088$ & $0.707$ & $41\,924$  & 10 & 32 & 32 &  -0.4919  \\ \hline %5
     \rowcolor{yellow!40}mixed & $Q_2^c Q_1^c Q_1^c$ & 0.2 & $0.044$ & $0.707$ & $165\,764$ & 4 & 36 & 31 & -0.0825  \\ \hline %6
      \rowcolor{yellow!50}mixed & $Q_2^c Q_1^c Q_1^c$ & 0.2 & $0.022$ & $0.707$ &  $659\,204$ & 8 & 50 & 53 & -0.0824 \\ \hline %7
       \rowcolor{yellow!60}mixed & $Q_2^c Q_1^c Q_1^c$ & 0.2 & $0.011$ & $0.707$ & $2\,629\,124$ & 8 & 79 & 38 & -0.0815 \\ \hline
  %8
  \hline
%\rowcolor{green!10}  $Q_1^c Q_1^{\text{stab}} Q_1^c$ & 0.2 & $0.353$ & $0.707$ & $1\,188$ & 5 & 2 & 3 &  -0.3731  \\ \hline %3
%\rowcolor{green!20}$Q_1^c Q_1^{\text{stab}} Q_1^c$ & 0.2 & $0.176$  & $0.707$ & $44\,204$ & 5 & 14 &3  &  -0.5127   \\ \hline %4
%\rowcolor{green!30}$Q_1^c Q_1^{\text{stab}} Q_1^c$& 0.2 & $0.088$ & $0.707$ &  &  &  &  &   \\ \hline %5
\hline
primal~\cite{heister2015primal} & $Q_2^c Q_1^c$ & 0.2 & $0.353$ & $0.707$ & $2\,507$ & 1 & - & 3 &  -0.3368  \\ \hline %3
primal~\cite{heister2015primal} &$Q_2^c Q_1^c$ & 0.2 & $0.176$  & $0.707$ & $9\,619$  & 1 & - & 3 &  -0.4479  \\ \hline %4
primal~\cite{heister2015primal} &$Q_2^c Q_1^c$& 0.2 & $0.088$ & $0.707$ & $37\,667$  & 5 &-  & 5 &   -0.4434 \\ \hline %5
primal~\cite{heister2015primal} &$Q_2^c Q_1^c$ & 0.2 & $0.044$ & $0.707$ & $149\,059$& 5 & - & 38 & -0.0818 \\ \hline %6
primal~\cite{heister2015primal} &$Q_2^c Q_1^c$ & 0.2 & $0.022$ & $0.707$ & $593\,027$ & 7 & - & 35 &  -0.0820  \\ \hline %7
primal~\cite{heister2015primal} &$Q_2^c Q_1^c$ & 0.2 & $0.011$ & $0.707$ & $2\,365\,699$ & 8 & - & 35 &  -0.0810  \\ \hline %8
  \hline
 \rowcolor{blue!10}mixed & $Q_2^c Q_1^c Q_1^c$ & 0.4999 & $0.353$ & $0.707$ & $2\,804$ & 10 & 24 &  3&   -0.2181 \\ \hline %3
  \rowcolor{blue!20}mixed & $Q_2^c Q_1^c Q_1^c$ & 0.4999 & $0.176$ & $0.707$ & $10\,724$  & 9 & 25 & 3 & -0.2869\\ \hline %4
    \rowcolor{blue!30}mixed & $Q_2^c Q_1^c Q_1^c$ & 0.4999& $0.088$ & $0.707$  & $41\,924$  & 6 & 32 & 29 & -0.1295  \\ \hline %5
  \rowcolor{blue!40}mixed & $Q_2^c Q_1^c Q_1^c$ & 0.4999 & $0.044$ & $0.707$ & $165\,764$ & 7 & 38 & 36 &  -0.0576 \\ \hline %6 
      \rowcolor{blue!50}mixed & $Q_2^c Q_1^c Q_1^c$ & 0.4999& $0.022$ & $0.707$ & $659\,204$  & 10 & 52 & 38 & -0.0585 \\ \hline %7
  \rowcolor{blue!60}mixed & $Q_2^c Q_1^c Q_1^c$& 0.4999 & $0.011$ & $0.707$ & $2\,629\,124$  & 11 & 80 & 41 &  -0.0578  \\ \hline %8 
  \hline
primal~\cite{heister2015primal} &$Q_2^c Q_1^c$ & 0.4999 & $0.353$ & $0.707$ & $2\,507$ & 1 & - & 3 &  -0.2077  \\ \hline %3
primal~\cite{heister2015primal} &$Q_2^c Q_1^c$ & 0.4999 & $0.176$  & $0.707$ & $9\,619$  & 1 & - & 3 &  -0.2788  \\ \hline %4
primal~\cite{heister2015primal} &$Q_2^c Q_1^c$& 0.4999 & $0.088$ & $0.707$ & $37\,667$  & 4 & - & 4 &  -0.2789  \\ \hline %5
primal~\cite{heister2015primal} &$Q_2^c Q_1^c$ & 0.4999 & $0.044$ & $0.707$ & $149\,059$& 5 & - &31  &  -0.0587 \\ \hline %6
primal~\cite{heister2015primal} &$Q_2^c Q_1^c$ & 0.4999 & $0.022$ & $0.707$ & $593\,027$ & 7 & - &  31&  -0.0583  \\ \hline %7
primal~\cite{heister2015primal} &$Q_2^c Q_1^c$ & 0.4999 & $0.011$ & $0.707$ & $2\,365\,699$ & 8& - & 34 & -0.5076  \\ \hline %8
\end{tabular}
\end{table}

  \begin{table}[htbp!]
  \scriptsize
\centering
%\textbf{Hanging block slit: robustness in $h$, $\lambda$ and $\epsilon$ for $\nu=0.5$; mixed}\\
%\vspace{0.2cm}
%\renewcommand*{\arraystretch}{1.2}
\caption{A hanging block with an initial slit for $\nu=0.5$, uniform mesh refinement, $Q_2^c Q_1^c Q_1^c$ elements, $\epsilon$ fixed and $\epsilon=2\,h$, $\kappa=10^{-2}$. Problem size in $\#$DoFs, average number of GMRES iterations ($\varnothing$lin) per Newton/active set (AS) step, the average number of inner CG iterations ($\varnothing$CG) per linear iteration, number of Newton/AS ($\#$AS), and goal functional displacement in a certain point ($u_y(0,1.99)$). Applied force $f=-8\cdot 10^{-7}\,\mathrm{N/mm^2}$.
}\label{HangingBlockSlitTable05}
\begin{tabular}{|l|l|c|c|c|r|r|r|r|l|} \hline
\multicolumn{10}{|c|}{Hanging block slit: robustness in $h$, $\lambda$ and $\epsilon$ for $\nu=0.5$; mixed; $\kappa=10^{-2}$}\\ \hline\hline
 \rowcolor{gray!50} \multicolumn{1}{|c|}{model} & \multicolumn{1}{|c|}{FE} & \multicolumn{1}{|c|}{$\nu$} & \multicolumn{1}{|c|}{$h$} & \multicolumn{1}{|c|}{$\epsilon$} & \multicolumn{1}{|c|}{$\#$DoFs} & \multicolumn{1}{|c|}{$\varnothing$lin} & \multicolumn{1}{|c|}{$\varnothing$CG} & \multicolumn{1}{|c|}{$\#$AS} & \multicolumn{1}{|c|}{$u_y(0,1.99)$}  \\
 \hline
 \rowcolor{Maroon!10}mixed & $Q_2^c Q_1^c Q_1^c$ & 0.5 & $0.353$ & $0.707$ & $2\,804$ & 9 & 23 & 3 &  -0.0578  \\ \hline %3
  \rowcolor{Maroon!20}mixed & $Q_2^c Q_1^c Q_1^c$ & 0.5 & $0.176$ & $0.707$ & $10\,724$  & 9 & 24 & 3 & -0.2835 \\ \hline %4
  \rowcolor{Maroon!30}mixed & $Q_2^c Q_1^c Q_1^c$ & 0.5 & $0.088$ & $0.707$ & $41\,924$  & 7 & 32 & 33 &  -0.0955 \\ \hline %5
  \rowcolor{Maroon!40}mixed & $Q_2^c Q_1^c Q_1^c$ & 0.5 & $0.044$ & $0.707$ & $165\,764$ & 6 & 37 & 38 &  -0.0584 \\ \hline %6
    \rowcolor{Maroon!50}mixed & $Q_2^c Q_1^c Q_1^c$ & 0.5 & $0.022$ & $0.707$ & $658\,436$  & 9 & 53 & 36 &  -0.0583  \\ \hline %7
  \rowcolor{Maroon!60}mixed & $Q_2^c Q_1^c Q_1^c$ & 0.5 & $0.011$ & $0.707$ & $2\,629\,124$ & 11 & 80 & 39 & -0.0579  \\ \hline %8
    \hline
 \rowcolor{magenta!10}mixed & $Q_2^c Q_1^c Q_1^c$ & 0.5 & $0.353$ & $0.707$ & $2\,804$ & 9 & 23 & 3 & -0.2166 \\ \hline %3
  \rowcolor{magenta!20}mixed & $Q_2^c Q_1^c Q_1^c$ & 0.5 & $0.176$ & $0.353$ & $10\,724$  & 7 & 25 & 4 &  -0.1033  \\ \hline %4
  \rowcolor{magenta!30}mixed & $Q_2^c Q_1^c Q_1^c$ & 0.5 & $0.088$ & $0.176$ & $41\,924$  & 6 & 30 & 14  &  -0.0701  \\ \hline %5
  \rowcolor{magenta!40}mixed & $Q_2^c Q_1^c Q_1^c$ & 0.5 & $0.044$ & $0.088$ & $165\,764$ & 5 & 36 & 109 &  -0.0572   \\ \hline %6
    \rowcolor{magenta!50}mixed & $Q_2^c Q_1^c Q_1^c$ & 0.5 & $0.022$ & $0.044$ & $658\,436$  & 7 & 40 & 805 &   -0.0516  \\ \hline %7
  %\rowcolor{magenta!60}mixed & $Q_2^c Q_1^c Q_1^c$ & 0.5 & $0.011$ & $0.022$ & $2\,629\,124$ &  &  &  &     \\ \hline %8
\end{tabular}
\end{table}

Tables~\ref{HangingBlockSlitTable} and~\ref{HangingBlockSlitTable05} show the iteration numbers of numerical tests for the hanging block with a slit for three Poisson ratios $\nu$ and $h$ refinement.
For the incompressible limit $\nu=0.5$, Table~\ref{HangingBlockSlitTable05} presents the results for $\epsilon$ fixed, and further in the pink rows, results for $\epsilon=2\,h$ are listed. The nearly constant number of GMRES iterations confirms the robustness in $\epsilon$ for $\nu=0.5$, tested for the hanging block with a slit on five levels of uniform refined meshes; see the last five rows in Table~\ref{HangingBlockSlitTable05}. 

\begin{remark}[High iteration numbers in the primal-dual active set method]
In Table~\ref{HangingBlockSlitTable05} in the pink rows, many active set/Newton iterations are required for $\epsilon\to 0$. Here, not the Poisson ratio is responsible, but the refinement in $h$ and $\epsilon$. For finer meshes with small $\epsilon$, the active set algorithm oscillates between a certain non-equal number of active nodes from the constraint. This effect leads to high total Newton iterations, even if the Newton algorithm converges fast; see also~\cite[Figure 14]{heister2015primal}.
\end{remark}
The number of CG iterations does not depend significantly on the size of $\kappa$ for this test setup. Further, the required CG iterations seem to be independent of $\lambda$ but sensitive to the mesh size. Aside from the robustness in $h$ and $\lambda$, we confirm the robustness in $\kappa$ for the hanging block test with a slit. Details on that can be found in~\cite[page 105]{Ma21_phd}.

%%%%%%%%%%%%%%%%%%%%%%%%%%%%%%%%%%%%%%%%%%%%%%%%%%%%%%%%%%%%%%%%%%%%%%%%%%%%%%%%%%%%%%%%%
%%%%%%%%%%%%%%%%%%%%%%%%%%%%%%%%%%%%%%%%%%%%%%%%%%%%%%%%%%%%%%%%%%%%%%%%%%%%%%%%%%%%%%%%%
%\newpage
\subsection{Sneddon's pressure-driven cavity}\label{sec_sneddon}

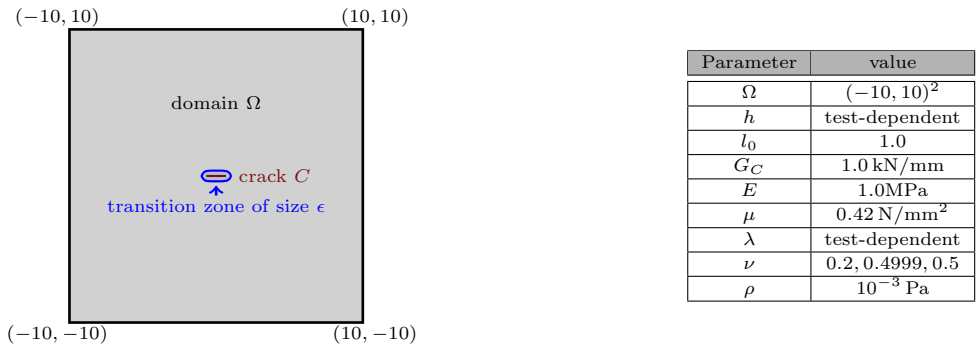
\begin{figure}[htbp!]
  \scriptsize %  \tiny
\begin{minipage}{0.6\textwidth}
\centering
\begin{tikzpicture}[xscale=0.65,yscale=0.65]
\draw[fill=gray!30] (-1,-1)  -- (-1,5) -- (5,5)  -- (5,-1) -- cycle;
\node  at (-1.25,5.25) {$(-10,10)$};
\node  at (-1.25,-1.25) {$(-10,-10)$};
\node at (5.25,-1.25) {$(10,-10)$};
\node at (5.25,5.25) {$(10,10)$};
\node at (2.0,3.5) {domain $\Omega$};
\draw [fill,opacity=0.1, blue] plot [smooth] coordinates { (1.7,2) (1.8,1.9) (2.2,1.9) (2.3,2) (2.2,2.1) (1.8,2.1) (1.7,2)};
\draw [thick, blue] plot [smooth] coordinates { (1.7,2) (1.8,1.9) (2.2,1.9) (2.3,2) (2.2,2.1) (1.8,2.1) (1.7,2)};
\draw[thick, draw=Maroon] (1.8,2)--(2.2,2);
\node[Maroon] at (3.2,2) {crack $C$};
\draw[->,blue] (2,1.6)--(2,1.8);
\node[blue] at (2,1.4) {transition zone of size $\epsilon$};
\end{tikzpicture}
%\caption{Geometry of the two-dimensional Sneddon's test in 2d. Two-dimensional domain of size $(-10,10)^2$ with an initial slit with a given pressure. }\label{SneddonGeo2d}
\end{minipage}
\hspace{0.2cm}
\begin{minipage}{0.35\textwidth}
\centering
%\begin{table}
\begin{tabular}{|c|c|}\hline
 \rowcolor{gray!50}\multicolumn{1}{|c|}{Parameter} &  \multicolumn{1}{c|}{value} \\ \hline\hline
$\Omega$ & $(-10,10)^2$  \\ \hline
$ h $  & test-dependent\\\hline
$ l_0 $ & $1.0$ \\ \hline
$ G_C $ & $1.0\,\mathrm{kN/mm}$ \\ \hline
$ E $ & $1.0\mathrm{MPa}$ \\ \hline
$\mu$ &   $0.42\,\mathrm{N/mm^2}$ \\ \hline
$\lambda$  & test-dependent \\ \hline
$ \nu $  & $0.2, 0.4999, 0.5$\\ \hline
$\rho$ & $10^{-3}\,\mathrm{Pa}$\\ \hline
\end{tabular}
%\end{table}
\end{minipage}
\caption{Left: geometry of the two-dimensional Sneddon's test in 2d. Domain of size $(-10,10)^2$ with a pressurized fracture. Right: setting of material and numerical parameters for Sneddon's benchmark test.}
\label{sneddon_params}
 \end{figure}

As a second example, we consider a benchmark test~\cite{schroder2021selection}, which is motivated by the book of Sneddon~\cite{sneddon1946distribution} and Sneddon and Lowengrub~\cite{SneddLow69}. We restrict ourselves to a 1d fracture $C$ on a 2d domain $\Omega = (-10,10)^2$ as depicted on the left in Figure~\ref{sneddon_params}. In this domain, an initial crack with length $2l_0 = 2.0$ and thickness $h$
 of two cells is prescribed with the help of the phase-field function $\varphi$, i.e.,
 $\varphi = 0$ in the crack and $\varphi = 1$ elsewhere.
 As boundary conditions, the displacements $u$ are set to zero on $\partial
 \Omega$.
 We use
homogeneous Neumann conditions for the phase-field variable, i.e., $\epsilon \partial_n \varphi = 0$ on $\partial \Omega$. 
The driving force is given by a constant pressure $\rho = 10^{-3}\,\mathrm{Pa}$ in the interior of the crack. An overview of the parameter setting is given in Figure~\ref{sneddon_params} on the right. 

 Two quantities of interest are discussed: the crack opening displacement (COD) and the total crack volume (TCV). 
The analytical solution (from~\cite{SneddLow69}) can be computed via 
\begin{align*}
\text{COD}_{\text{ref}} = 2\frac{pl_0}{E'}\left(1-\frac{x^2}{l_0^2} \right)^\frac{1}{2},
\end{align*}
where $E' \coloneqq  \frac{E}{1-\nu^2}$, $E$ is the Young modulus and $\nu$ is the Poisson ratio.
The TCV can be computed numerically with 
\begin{align*}
\text{TCV} = \int_\Omega u(x,y) \cdot \nabla \varphi (x,y) d(x,y).
\end{align*}
The analytical solution (from~\cite{SneddLow69}) is given by 
\begin{align*}
\text{TCV}_{\text{ref}} = \frac{2\pi p l_0^2}{E'}.
\end{align*}
%In Table~\ref{reference_cod_tcv}, the manufactured reference values (on an infinite domain) of COD and TCV are given for the three considered Poisson ratios.\\
%The idea of a consistently stable equal-order discretization is repeated briefly in the following.

In Table~\ref{Sneddon2dKappSmall}, for $\kappa=10^{-8}$, the average number of CG iterations increases with a decreasing
mesh size. We observe an increase in the CG iteration numbers in particular for the incompressible limit $\nu=0.5$ and finer meshes, where we finally do not get convergence in the solver for smaller $h$. Already for $\nu=0.4999$ and a problem size of less than $300\,000$ DoFs, the average number of CG iterations is above $100$. 
%%%%%%%%%%%%%%%%%%%%%%%%%%%%%%%%%%%%%
\begin{remark}[Difficulties considering small $\kappa$]\label{kappa_problems}
In Table~\ref{Sneddon2dKappSmall}, compared to Table~\ref{Sneddon2d}, we can evaluate the impact of the setting of $\kappa$. We compute Sneddon's test for different mesh sizes $h$, fixed bandwidth $\epsilon$, for three Poisson ratios $\nu=0.2, 0.4999$, and $\nu=0.5$, and for a small and large regularization parameter $\kappa=10^{-2}$ and $\kappa=10^{-8}$ to evaluate its impact on the behavior of the CG solver. These solver dependencies on $\kappa$ have a natural correspondence in error estimates. For a decoupled
linearized system, such estimates are shown in~\cite[Section 5.5]{wick2020multiphysics}. A numerical error analysis for this test on a good choice of $\kappa$ can be found in~\cite{KOLDITZ2022100047}.

Further, we observe an increased number of CG iterations for high Poisson ratios. The number of GMRES and AS iterations do not differ significantly for different $\kappa$.
\end{remark}

  \begin{table}[htbp!]
  \scriptsize %  \tiny
\centering
\caption{Sneddon's pressure-driven cavity in 2d. Average number of GMRES iterations ($\#$lin) per Newton step ($\#$AS), the average number of CG iterations ($\varnothing$CG) per linear iteration. Based on the newly developed mixed model with $Q_2^c Q_1^c Q_1^c$ elements for different problem sizes and setting of the length scale parameter $\epsilon$ for three Poisson ratios. Quantities of interest: COD$_{\text{max}}$ and TCV and $\kappa =10^{-8}$. Uniform refined meshes.
}\label{Sneddon2dKappSmall}
\begin{tabular}{|l|c|c|c|r|r|r|r|l|r|} \hline
\multicolumn{10}{|c|}{Sneddon's test: robustness in $h$, $\lambda$, $\epsilon$; $\kappa=10^{-8}$; mixed}\\ \hline\hline
 \rowcolor{gray!50} \multicolumn{1}{|c|}{FE} & \multicolumn{1}{|c|}{$\nu$} & \multicolumn{1}{|c|}{$h$} & \multicolumn{1}{|c|}{$\epsilon$} & \multicolumn{1}{|c|}{$\#$DoFs} & \multicolumn{1}{|c|}{$\varnothing$lin} & \multicolumn{1}{|c|}{$\varnothing$CG} & \multicolumn{1}{|c|}{$\#$AS} & \multicolumn{1}{|c|}{COD$_{\max}$} & \multicolumn{1}{|c|}{TCV} \\
 \hline
  \hline
   \rowcolor{yellow!10}$Q_2^c Q_1^c Q_1^c$ & 0.2 & $0.707$ & $1.414$ & $16\,484$  & 3 & 26 & 4 & 0.00282 & 0.0240\\ \hline %2
  \rowcolor{yellow!20}$Q_2^c Q_1^c Q_1^c$ & 0.2 & $0.353$  & $1.414$ & $64\,964$   & 6 & 28 & 6 & 0.00270 & 0.0189\\ \hline %3
    \rowcolor{yellow!30}$Q_2^c Q_1^c Q_1^c$ & 0.2 & $0.176$ & $1.414$ & $257\,924$ & 9 & 35 & 4 & 0.00260 & 0.0164\\ \hline %4
  \rowcolor{yellow!40}$Q_2^c Q_1^c Q_1^c$ & 0.2 & $0.088$ & $1.414$ & $1\,027\,844$ & 12 & 31 & 5 & 0.00252 & 0.0150\\ \hline %5
 % \hline
%  \rowcolor{green!10}$Q_1^c Q_1^{\text{stab}} Q_1^c$ & 0.2 & $0.707$ & $1.414$ & $6\,724$ & 6 & 8 & 2 & 0.00280 & 0.0235\\ \hline %2
%   \rowcolor{green!20}$Q_1^c Q_1^{\text{stab}} Q_1^c$ & 0.2 & $0.353$  & $1.414$ & $26\,244$ & 8 & 9 & 3 & 0.00269 & 0.0188\\ \hline %3
%     \rowcolor{green!30}$Q_1^c Q_1^{\text{stab}} Q_1^c$ & 0.2 & $0.176$ & $1.414$ & $103\,684$ & 9 & 11 & 3 & 0.00260 & 0.0164\\ \hline %4
%   \rowcolor{green!40}$Q_1^c Q_1^{\text{stab}} Q_1^c$ & 0.2 & $0.088$ & $1.414$ & $412\,164$ & 9 & 12 & 5 & 0.00252 & 0.0150\\ \hline %5
   ref.~\cite{SneddLow69} & 0.2 &  &  & & &  & & 0.0019200 & 0.00603\\ \hline
  \hline
 \rowcolor{blue!10}$Q_2^c Q_1^c Q_1^c$ & 0.4999 & $0.707$ & $1.414$ & $16\,484$  & 3 & 31 & 6 & 3.0383e-05 & 0.000257\\ \hline %2
  \rowcolor{blue!20}$Q_2^c Q_1^c Q_1^c$ & 0.4999 & $0.353$ & $1.414$ & $64\,964$  & 7 & 46 & 8 & 3.6024e-05 & 0.000254\\ \hline %3
    \rowcolor{blue!30}$Q_2^c Q_1^c Q_1^c$ & 0.4999& $0.176$ & $1.414$ & $257\,924$ & 6 & 107 & 39 &  3.9899e-05 & 0.000252 \\ \hline %4
  \rowcolor{blue!40}$Q_2^c Q_1^c Q_1^c$ & 0.4999 & $0.088$ & $1.414$ & $1\,027\,844$  & 5 & 57 & 24 & 4.2265e-05 & 0.000250\\ \hline %5 
            ref.~\cite{SneddLow69} & 0.4999 &  &  & & &  &  & 0.0015001 & 0.004713\\ \hline 
        \hline
 \rowcolor{Maroon!10}$Q_2^c Q_1^c Q_1^c$ & 0.5 & $0.707$ & $1.414$ & $16\,484$  & 3 & 31 & 3 & 2.9937e-20 & 7.1504e-20 \\ \hline %2
  \rowcolor{Maroon!20}$Q_2^c Q_1^c Q_1^c$ & 0.5 & $0.353$ & $1.414$ & $64\,964$  & 6 & 25 & 2 & 1.3258e-19 &2.3835e-19 \\ \hline %3
  \rowcolor{Maroon!30}$Q_2^c Q_1^c Q_1^c$ & 0.5 & $0.176$ & $1.414$ & $257\,924$  & 5 & 59 & 7 & 1.9309e-19 & 7.8981e-19\\ \hline %4
  \rowcolor{Maroon!40}$Q_2^c Q_1^c Q_1^c$ & 0.5 & $0.088$ & $1.414$ & $1\,027\,844$  & -  & - & - &- & - \\ \hline %5 keine Konvergenz
    \hline
 \rowcolor{magenta!10}$Q_2^c Q_1^c Q_1^c$ & 0.5 & $0.707$ & $1.414$ & $16\,484$  & 11 & 37 & 3 &2.4585e-15  &  1.2562e-14\\ \hline %2
  \rowcolor{magenta!20}$Q_2^c Q_1^c Q_1^c$ & 0.5 & $0.353$ & $0.707$ & $64\,964$  & 6 & 32 & 3 & 2.3632e-18 & 1.0069e-17 \\ \hline %3
  \rowcolor{magenta!30}$Q_2^c Q_1^c Q_1^c$ & 0.5 & $0.176$ & $0.353$ & $257\,924$ & 10 & 30 & 3 & 6.5953e-18 & 1.4749e-16\\ \hline %4
  \rowcolor{magenta!40}$Q_2^c Q_1^c Q_1^c$ & 0.5 & $0.088$ & $0.176$ & $1\,027\,844$   & 14 & 38 & 3 & 1.2397e-18 & 2.6778e-18\\ \hline %5
  ref.~\cite{SneddLow69} & 0.5 &  &  & & & &  &  0.0015000 & 0.0047124\\ \hline
\end{tabular}
\end{table}

 This observation is confirmed by the numerical results from Table~\ref{Sneddon2dKappSmallCgforS}, where a CG solver preconditioned with AMG is used to approximate $\hat{S}^{-1}$. The numerical results in Table~\ref{Sneddon2dKappSmallCgforS} are based on the same tests as in Table~\ref{Sneddon2dKappSmall} but for $\nu=0.4999$ and $\nu=0.5$. The number of linear iterations is moderate, and at most six CG iterations are needed for $\hat{S}^{-1}$. 

  \begin{table}[htbp!]
  \scriptsize %  \tiny
  \caption{Sneddon's pressure-driven cavity in 2d. Average number of GMRES iterations ($\#$lin) per Newton step ($\#$AS), the average number of CG iterations ($\varnothing$CG) per linear iteration, CG plus AMG is used for $(g(\tilde{\varphi})A_u)^{-1}$ and $\hat{S}^{-1}$. Based on the newly developed mixed model with $Q_2^c Q_1^c Q_1^c$ elements for different problem sizes and setting of the length scale parameter $\epsilon$ for two Poisson ratios. Uniform refined meshes. 
}\label{Sneddon2dKappSmallCgforS}
\centering
\begin{tabular}{|l|c|c|c|r|r|r|c|l|} \hline
\multicolumn{9}{|c|}{Sneddon's test: robustness in $h$, $\lambda$, $\epsilon$; $\kappa=10^{-8}$; mixed; CG+AMG for $\hat{S}^{-1}$}\\ \hline\hline
 \rowcolor{gray!50} \multicolumn{1}{|c|}{FE} & \multicolumn{1}{|c|}{$\nu$} & \multicolumn{1}{|c|}{$h$} & \multicolumn{1}{|c|}{$\epsilon$} & \multicolumn{1}{|c|}{$\#$DoFs} & \multicolumn{1}{|c|}{$\varnothing$lin} & \multicolumn{1}{|c|}{$\varnothing$CG $(g(\tilde{\varphi})A_u)^{-1}$} & \multicolumn{1}{|c|}{$\varnothing$CG $\hat{S}^{-1}$} & \multicolumn{1}{|c|}{$\#$AS}\\
 \hline
  \hline
 \rowcolor{blue!10}$Q_2^c Q_1^c Q_1^c$ & 0.4999 & $0.707$ & $1.414$ & $16\,484$  & 3 & 26 & 1 & 3 \\ \hline %2
  \rowcolor{blue!20}$Q_2^c Q_1^c Q_1^c$ & 0.4999 & $0.353$ & $1.414$ & $64\,964$ & 8 & 56 & 6 & 8 \\ \hline %3
    \rowcolor{blue!30}$Q_2^c Q_1^c Q_1^c$ & 0.4999& $0.176$ & $1.414$ & $257\,924$ & 6 & 106 & 6 & 38 \\ \hline %4
  \rowcolor{blue!40}$Q_2^c Q_1^c Q_1^c$ & 0.4999 & $0.088$ & $1.414$ & $1\,027\,844$ & 6 & 42 & 6 & 69 \\ \hline %5 
       %     ref.~\cite{SneddLow69} & 0.4999 &  &  & & &  &  & 0.0015001 & 0.004713\\ \hline 
        \hline
 \rowcolor{red!10}$Q_2^c Q_1^c Q_1^c$ & 0.5 & $0.707$ & $1.414$ & $16\,484$ & 10 & 36 & 1 & 3 \\ \hline %2
  \rowcolor{red!20}$Q_2^c Q_1^c Q_1^c$ & 0.5 & $0.353$ & $1.414$ & $64\,964$ & 6 & 26 & 6 & 8\\ \hline %3
  \rowcolor{red!30}$Q_2^c Q_1^c Q_1^c$ & 0.5 & $0.176$ & $1.414$ & $257\,924$ & 6 & 63 & 6 & 37\\ \hline %4
  \rowcolor{red!40}$Q_2^c Q_1^c Q_1^c$ & 0.5 & $0.088$ & $1.414$ & $1\,027\,844$ & 7 & 41 & 6 & 101\\ \hline %5 
   % \hline
%% \rowcolor{magenta!10}$Q_2^c Q_1^c Q_1^c$ & 0.5 & $0.707$ & $1.414$ & $16\,484$ & 11 &  &  & 3\\ \hline %2
 % \rowcolor{magenta!20}$Q_2^c Q_1^c Q_1^c$ & 0.5 & $0.353$ & $0.707$ & $64\,964$ & 7 &  &  & 3\\ \hline %3
 % \rowcolor{magenta!30}$Q_2^c Q_1^c Q_1^c$ & 0.5 & $0.176$ & $0.353$ & $257\,924$ & 8 &  &  & 3\\ \hline %4
 % \rowcolor{magenta!40}$Q_2^c Q_1^c Q_1^c$ & 0.5 & $0.088$ & $0.176$ & $1\,027\,844$ & 13 &  &  & 3\\ \hline %5
%  ref.~\cite{SneddLow69} & 0.5 &  &  & & & &  &  0.0015000 & 0.0047124 & \\ \hline
\end{tabular}
\end{table}

%%%%%%%%%%%%%%%%%%%%%%%%%%%%%%%%%%%%%
As expected in Tables~\ref{Sneddon2dKappSmall},~\ref{Sneddon2dKappSmallCgforS} and~\ref{Sneddon2d}, considering the quantities of interest COD$_{\text{max}}$ and TCV, they get vanishingly small for high Poisson ratios. This is what we expected for incompressible solids: a closed domain does not change its volume; the opening of the initial crack in the interior of the domain is avoided. For $\nu=0.2$, the quantities of interest are acceptable compared to the reference values. % from Table~\ref{reference_cod_tcv}. 
Also for $\nu=0.2$, since all computations are conducted with uniformly refined meshes, moderate problem sizes, and fixed $\epsilon$, we cannot expect excellent results in the quantities of interest. 

%%%%%%%%%%%%%%%%%%%%%%%%%%%%%%%%%%%%%%%%%%%%%%%%%%%%%%%%%%%%%%%%%%%%%%%%

  \begin{table}[htbp!]
  \scriptsize %  \tiny
\centering
\caption{Sneddon's pressure-driven cavity. Average number of GMRES iterations ($\#$lin) per Newton step ($\#$AS), average number of CG iterations ($\varnothing$CG) per linear iteration. Based on the newly developed mixed model with $Q_2^c Q_1^c Q_1^c$ elements for different problem sizes and setting of the length scale parameter $\epsilon$ for three Poisson ratios. Quantities of interest: COD$_{\text{max}}$ and TCV and $\kappa =10^{-2}$. Uniform refined meshes.
}\label{Sneddon2d}
\begin{tabular}{|l|c|c|c|r|r|r|r|l|r|} \hline
\multicolumn{10}{|c|}{Sneddon's test: robustness in $h$, $\lambda$, $\epsilon$; $\kappa=10^{-2}$; mixed}\\ \hline\hline
 \rowcolor{gray!50} \multicolumn{1}{|c|}{FE} & \multicolumn{1}{|c|}{$\nu$} & \multicolumn{1}{|c|}{$h$} & \multicolumn{1}{|c|}{$\epsilon$} & \multicolumn{1}{|c|}{$\#$DoFs} & \multicolumn{1}{|c|}{$\varnothing$lin} & \multicolumn{1}{|c|}{$\varnothing$CG} & \multicolumn{1}{|c|}{$\#$AS} & \multicolumn{1}{|c|}{COD$_{\max}$} & \multicolumn{1}{|c|}{TCV} \\
 \hline
  \hline
   \rowcolor{yellow!10}$Q_2^c Q_1^c Q_1^c$ & 0.2 & $0.707$ & $1.414$ & $16\,484$  & 2 & 16 & 4 & 0.00248 & 0.0224\\ \hline %2 % 1.31345485e-05
  \rowcolor{yellow!20}$Q_2^c Q_1^c Q_1^c$ & 0.2 & $0.353$  & $1.414$ & $64\,964$   & 8 & 18 & 4 & 0.00227 & 0.0173 \\ \hline %3 % 1.04672475e-05
    \rowcolor{yellow!30}$Q_2^c Q_1^c Q_1^c$ & 0.2 & $0.176$ & $1.414$ & $257\,924$ & 9 & 18 & 15 & 0.00206 & 0.0145\\ \hline %4 % 9.02498973e-06
  \rowcolor{yellow!40}$Q_2^c Q_1^c Q_1^c$ & 0.2 & $0.088$ & $1.414$ & $1\,027\,844$ & 15 & 28 & 5 & 0.00190 & 0.0129\\ \hline %5 % 8.21127794e-06
 % \hline
%  \rowcolor{green!10}$Q_1^c Q_1^{\text{stab}} Q_1^c$ & 0.2 & $0.707$ & $1.414$ & $6\,724$ & 6 & 8 & 2 & 0.00245 & 0.0221 \\ \hline %2 % 1.30194625e-05
%   \rowcolor{green!20}$Q_1^c Q_1^{\text{stab}} Q_1^c$ & 0.2 & $0.353$  & $1.414$ & $26\,244$ & 8 & 9 & 2 & 0.00224 & 0.0171\\ \hline %3 % 1.04234330e-05
%     \rowcolor{green!30}$Q_1^c Q_1^{\text{stab}} Q_1^c$ & 0.2 & $0.176$ & $1.414$ & $103\,684$ & 9 & 10 & 6 & 0.00205 & 0.0144\\ \hline %4 % 9.00953465e-06
%   \rowcolor{green!40}$Q_1^c Q_1^{\text{stab}} Q_1^c$ & 0.2 & $0.088$ & $1.414$ & $412\,164$ & 9 & 12 & 8 & 0.00190 & 0.0129 \\ \hline %5 8.20615293e-06
   ref.~\cite{SneddLow69} & 0.2 &  & & &  &  & & 0.0019200 & 0.0060\\ \hline
  \hline
 \rowcolor{blue!10}$Q_2^c Q_1^c Q_1^c$ & 0.4999 & $0.707$ & $1.414$ & $16\,484$  & 13 & 16 & 3 & 3.0833e-05 & 0.000269\\ \hline %2 % 4.11989195e-06
  \rowcolor{blue!20}$Q_2^c Q_1^c Q_1^c$ & 0.4999 & $0.353$ & $1.414$ & $64\,964$  & 8 & 18 & 14 & 3.1739e-05 & 0.000242\\ \hline %3 % 2.86516721e-06
    \rowcolor{blue!30}$Q_2^c Q_1^c Q_1^c$ & 0.4999& $0.176$ & $1.414$ & $257\,924$ & 6 & 18 & 93 &  3.3667e-05 & 0.000224 \\ \hline %4 % 1.26622494e-06
  \rowcolor{blue!40}$Q_2^c Q_1^c Q_1^c$ & 0.4999 & $0.088$ & $1.414$ & $1\,027\,844$ & 7 & 26 & 65 & 3.4560e-05 & 0.000216 \\ \hline %5 % 9.75339950e-07
          ref.~\cite{SneddLow69} & 0.4999 &  &  & & &  &  & 0.0015001 & 0.004713\\ \hline 
        \hline
 \rowcolor{Maroon!10}$Q_2^c Q_1^c Q_1^c$ & 0.5 & $0.707$ & $1.414$ & $16\,484$  & 9 & 14 & 3 & 1.7339e-19 & 5.8895e-19\\ \hline %2 changed
  \rowcolor{Maroon!20}$Q_2^c Q_1^c Q_1^c$ & 0.5 & $0.353$ & $1.414$ & $64\,964$  & 9 & 18 & 14 & 2.3734e-19 & 5.6268e-18\\ \hline %3 changed
  \rowcolor{Maroon!30}$Q_2^c Q_1^c Q_1^c$ & 0.5 & $0.176$ & $1.414$ & $257\,924$  & 11 & 18 & 14 & 5.6547e-20 & 6.0823e-18\\ \hline %4 changed
  \rowcolor{Maroon!40}$Q_2^c Q_1^c Q_1^c$ & 0.5 & $0.088$ & $1.414$ & $1\,027\,844$  & 5 & 26 & 39 & 7.7351e-19 & 2.2733e-17 \\ \hline %5
    \hline
 \rowcolor{magenta!10}$Q_2^c Q_1^c Q_1^c$ & 0.5 & $0.707$ & $1.414$ & $16\,484$  & 9 & 14 & 3 & 1.5881e-19 & 5.8895e-19\\ \hline %2 infty
  \rowcolor{magenta!20}$Q_2^c Q_1^c Q_1^c$ & 0.5 & $0.353$ & $0.707$ & $64\,964$  & 6 & 17 & 3 & 1.9290e-19 & 1.8057e-18\\ \hline %3 infty
  \rowcolor{magenta!30}$Q_2^c Q_1^c Q_1^c$ & 0.5 & $0.176$ & $0.353$ & $257\,924$  & 6 & 18 & 3 & 4.1847e-19 & 2.1156e-18\\ \hline %4 infty
  \rowcolor{magenta!40}$Q_2^c Q_1^c Q_1^c$ & 0.5 & $0.088$ & $0.176$ & $1\,027\,844$  & 10 & 26 & 3 & 2.4801e-18 &  9.5514e-18\\ \hline %5 infty
         ref.~\cite{SneddLow69} & 0.5 &  &  & & & &  &  0.0015000 & 0.0047124\\ \hline
\end{tabular}
\end{table}
In Table~\ref{Sneddon2d}, the same computations are conducted as in Table~\ref{Sneddon2dKappSmall} and Table~\ref{Sneddon2dKappSmallCgforS} for $\kappa=10^{-2}$ to discuss the statement of Remark~\ref{kappa_problems}.
The COD values are close to the reference values. Here, a large regularization parameter $\kappa=10^{-2}$ stabilizes the block $(g(\tilde{\varphi})A_u)^{-1}$. Further, the linear iterations are stable, and also the inner CG iterations are relatively constant. In the last four rows of Table~\ref{Sneddon2d}, similar to Table~\ref{HangingBlockSlitTable05}, results of four tests with $\epsilon=2\,h$ are listed to check the robustness in $\epsilon$ for $\nu=0.5$, which can be confirmed for Sneddon's benchmark test.

%%%%%%%%%%%%%%%%%%%%%%%%%%%%%%%%%%%%%%%%%%%%%%%%%%%%%%%%%%%%%%%%%%%%%%%%%%%%%%%%%%%%%%%%%
 \subsection{Sneddon's pressure-driven cavity, layered}\label{sec_sneddon_layered}

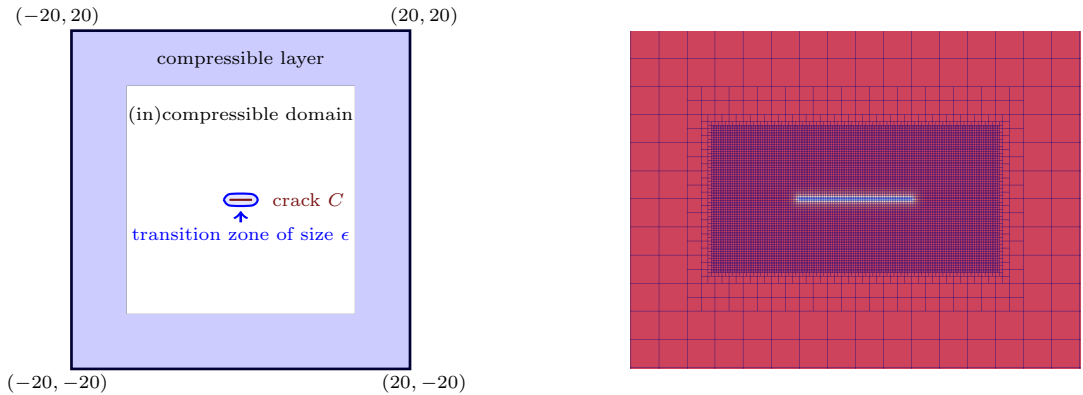
\begin{figure}[htbp!]
\scriptsize
\begin{minipage}{0.58\textwidth}
\centering
%\small
\begin{tikzpicture}[xscale=0.75,yscale=0.75]
\draw (-1,-1)  -- (-1,5) -- (5,5)  -- (5,-1) -- cycle;
\node  at (-1.25,5.25) {$(-20,20)$};
\node  at (-1.25,-1.25) {$(-20,-20)$};
\node at (5.25,-1.25) {$(20,-20)$};
\node at (5.25,5.25) {$(20,20)$};
\draw[gray, fill,opacity=0.2, blue] (-1,-1)  -- (-1,5) -- (5,5)  -- (5,-1) -- cycle;
\draw[fill=gray!30] (0,0)  -- (0,4) -- (4,4)  -- (4,0) -- cycle;
\draw[fill, white] (0,0)  -- (0,4) -- (4,4)  -- (4,0) -- cycle;
\node at (2.0,4.5) {compressible layer};
\node at (2.0,3.5) {(in)compressible domain};
\draw [fill,opacity=0.1, blue] plot [smooth] coordinates { (1.7,2) (1.8,1.9) (2.2,1.9) (2.3,2) (2.2,2.1) (1.8,2.1) (1.7,2)};
\draw [thick, blue] plot [smooth] coordinates { (1.7,2) (1.8,1.9) (2.2,1.9) (2.3,2) (2.2,2.1) (1.8,2.1) (1.7,2)};
\draw[thick, draw=Maroon] (1.8,2)--(2.2,2);
\node[Maroon] at (3.2,2) {crack $C$};
\draw[->,blue] (2,1.6)--(2,1.8);
\node[blue] at (2,1.4) {transition zone of size $\epsilon$};
\end{tikzpicture}
\end{minipage}
\hspace{0.2cm}
\begin{minipage}{0.37\textwidth}
% \end{figure}
% \begin{figure}
    \centering
    \includegraphics[width=0.98\textwidth]{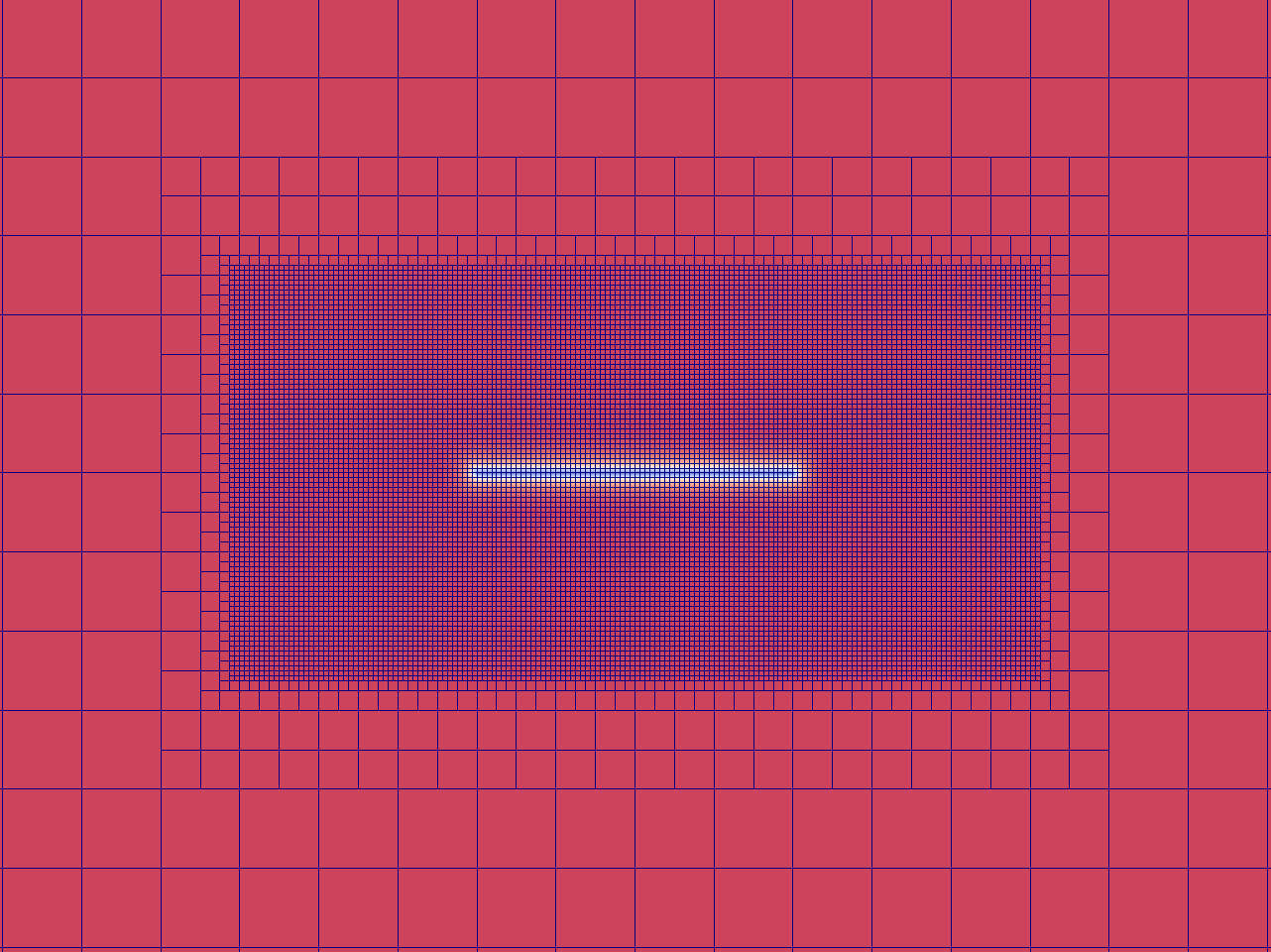}
    \end{minipage}
     %   \caption{Zoom-in to the geometrically refined mesh used in Table~\ref{Sneddon2dlayeredAdaptive}.}
    %\label{geo-adapt}
        \caption{Left: Geometry of the two-dimensional Sneddon's test with a compressible layer of size $10$. Further, the inside of the initial crack is assumed to be compressible~\cite{basava2020adaptive}. Right: Close zoom-in to the geometrically refined mesh around the crack, used in Table~\ref{Sneddon2dlayeredAdaptive}.}\label{SneddonGeoLayered}
\end{figure}

 As a fourth test case, the pressure-driven cavity from~\cite{schroder2021selection} is modified similarly to~\cite{basava2020adaptive}. We consider a two-dimensional domain $\Omega = (-20,20)^2$.
 %as sketched in Figure~\ref{SneddonGeoLayered}. 
 In contrast to the previous Sneddon test, a compressible layer of size $10$ is added around the incompressible domain to allow deforming of the solid on a finite domain. So the Poisson ratio changes over the domain for the layered Sneddon test. We expect to get better results concerning COD and TCV on a finite domain compared to the reference values on an infinite domain.
 A sketch of the geometry is given in Figure~\ref{SneddonGeoLayered} on the left. The setting of the material and numerical parameters is the same as in the previous section.

In Figure~\ref{SneddonGeoLayered} on the right, a zoom-in snapshot of the inner domain is given to show the geometric refinement for the tests in Tables~\ref{Sneddon2dlayered},~\ref{Sneddon2dlayeredAdaptive}, and~\ref{Sneddon2dlayeredAdaptiveCgforS}. Aside from the adaptively refined mesh, we set $\epsilon=h$, depending on the current mesh size. The total numbers of degrees of freedom ($\#$DoFs) on $\Omega$ are listed in the numerical results in Tables~\ref{Sneddon2dlayered} to~\ref{Sneddon2dlayeredAdaptiveCgforS}.

  \begin{table}[htbp!]
  \scriptsize %  \tiny
\centering
\caption{Sneddon's pressure-driven cavity layered. Average number of GMRES iterations ($\varnothing$lin) per Newton step ($\#$AS). Computations based on the newly developed mixed model with $Q_2^c Q_1^c Q_1^c$ elements for different problem size, $\epsilon=h$ for three Poisson ratios. Quantities of interest: COD$_{\text{max}}$ and TCV and $\kappa =10^{-2}$. Geometrically refined mesh in the area around the crack zone as depicted in Figure~\ref{SneddonGeoLayered}.
}\label{Sneddon2dlayered}
\begin{tabular}{|l|l|l|r|r|r|l|r|} \hline
\multicolumn{8}{|c|}{Sneddon layered adaptive: robustness in $h$, $\lambda$, $\epsilon$; $\epsilon=h$ $\kappa=10^{-2}$; mixed}\\ \hline\hline
 \rowcolor{gray!50} \multicolumn{1}{|c|}{FE} & \multicolumn{1}{|c|}{$\nu$} & \multicolumn{1}{|c|}{$h$} & \multicolumn{1}{|c|}{$\#$DoFs} & \multicolumn{1}{|c|}{$\varnothing$lin} & \multicolumn{1}{|c|}{$\#$AS} & \multicolumn{1}{|c|}{COD$_{\max}$} & \multicolumn{1}{|c|}{TCV} \\
 \hline
  \hline  % bisher gut
    \rowcolor{yellow!5}$Q_2^c Q_1^c Q_1^c$ & 0.2 & $0.353$ & $257\,924$   & 18  & 3 & 0.00214077 & 0.0097207 \\ \hline %4+1
  \rowcolor{yellow!10}$Q_2^c Q_1^c Q_1^c$ & 0.2 & $0.176$  & $263\,604$   & 25 & 4 & 0.00188194 & 0.0069353\\ \hline %4+2
  \rowcolor{yellow!15}$Q_2^c Q_1^c Q_1^c$ & 0.2 & $0.088$  & $282\,484$  & 20 & 3 & 0.00163459 & 0.0055415\\ \hline %4+3
  \rowcolor{yellow!20}$Q_2^c Q_1^c Q_1^c$ & 0.2 & $0.044$ & $350\,804$   & 17 & 3 & 0.00136002 & 0.0044906\\ \hline %4+4
    \rowcolor{yellow!25}$Q_2^c Q_1^c Q_1^c$ & 0.2 & $0.022$ & $610\,164$  & 19 & 4 & 0.00104379 & 0.0034504\\ \hline %4+5
        \rowcolor{yellow!30}$Q_2^c Q_1^c Q_1^c$ & 0.2 & $0.011$ & $1\,620\,244$  & 24 & 6 & 0.00071731 & 0.0024168\\ \hline %4+6
            \rowcolor{yellow!35}$Q_2^c Q_1^c Q_1^c$ & 0.2 & $0.0055$ & $5\,606\,324$  & 28 & 6 & 0.00043963 & 0.0015294\\ \hline %3+7
 ref.~\cite{SneddLow69} & 0.2 &  &  &  & & 0.00192000 & 0.0060318\\ \hline
   \hline  % bisher gut
           \rowcolor{blue!5}$Q_2^c Q_1^c Q_1^c$ & 0.4999 & $0.353$ & $257\,924$    &40  & 2 & 0.00205349 & 0.0108334\\ \hline %4+1
  \rowcolor{blue!10}$Q_2^c Q_1^c Q_1^c$ & 0.4999 & $0.176$  & $263\,604$    & 52 & 2 & 0.00168136 & 0.0069892\\ \hline %4+2
  \rowcolor{blue!15}$Q_2^c Q_1^c Q_1^c$ & 0.4999 & $0.088$  & $282\,484$    & 58  & 3 & 0.00143863 & 0.0052347\\ \hline %4+3
  \rowcolor{blue!20}$Q_2^c Q_1^c Q_1^c$ & 0.4999 & $0.044$ & $350\,804$    & 58 & 3 & 0.00122931 & 0.0041947\\ \hline %4+4
    \rowcolor{blue!25}$Q_2^c Q_1^c Q_1^c$ & 0.4999 & $0.022$ & $610\,164$   & 62 & 4 & 0.00099810 & 0.0033284\\ \hline %4+5
        \rowcolor{blue!30}$Q_2^c Q_1^c Q_1^c$ & 0.4999 & $0.011$ & $1\,620\,244$  & 150 & 4 & 0.00073786 & 0.0024681\\ \hline %4+6
        \rowcolor{blue!35}$Q_2^c Q_1^c Q_1^c$ & 0.4999 & $0.0055$ & $5\,606\,324$  & 318 & 7 & 0.00048545 & 0.0016595\\ \hline %4+7 
        ref.~\cite{SneddLow69} & 0.4999 &  &  &  &  & 0.00150019 & 0.0047130\\ \hline 
        \hline
        %%%%%%%%%%%%%%%%%%%%%%%%%%%%%%%%%%%%%%%%%%%%%  % bisher gut
  \rowcolor{Maroon!5}$Q_2^c Q_1^c Q_1^c$ & 0.5 & $0.353$ & $257\,924$  & 40 & 2 & 0.00205332 & 0.0108334\\ \hline %4+1
  \rowcolor{Maroon!10}$Q_2^c Q_1^c Q_1^c$ & 0.5 & $0.176$  & $263\,604$  & 52 & 2 & 0.00168117 & 0.0069889 \\ \hline %4+2
  \rowcolor{Maroon!15}$Q_2^c Q_1^c Q_1^c$ & 0.5 & $0.088$  & $282\,484$  & 59 & 3 & 0.00143847 & 0.0052343\\ \hline %4+3
    \rowcolor{Maroon!20}$Q_2^c Q_1^c Q_1^c$ & 0.5 & $0.044$ & $350\,804$    & 57 & 4 & 0.00122920 & 0.0041944\\ \hline %4+4
    \rowcolor{Maroon!25}$Q_2^c Q_1^c Q_1^c$ & 0.5 & $0.022$ & $610\,164$   & 65 & 4 & 0.00099804 & 0.0033282\\ \hline %4+5
        \rowcolor{Maroon!30}$Q_2^c Q_1^c Q_1^c$ & 0.5 & $0.011$ & $1\,620\,244$  & 155 & 4 & 0.00073784 & 0.0024681\\ \hline %4+6
        \rowcolor{Maroon!35}$Q_2^c Q_1^c Q_1^c$ & 0.5 & $0.0055$ & $5\,606\,324$  & 271  & 7 & 0.00048545 & 0.0016595\\ \hline %4+7 
       ref.~\cite{SneddLow69} & 0.5 &  &  &  &  &  0.00150000 & 0.0047124\\ \hline 
\end{tabular} % laeuft nochmal mit eps =h und 8 verfeinerungsschritten, hier auch cod plotten?
\end{table}

In Table~\ref{Sneddon2dlayered}, the results for the Sneddon test in 2d with a compressible layer around a possibly incompressible domain are given for three Poisson ratios and adaptively refined meshes, with $\epsilon=h$, and $\kappa=10^{-2}$. We choose $\kappa=10^{-2}$ to avoid the effects of $\kappa$ on the inner CG iterations. For large $\kappa$, the computed quantities of interest COD$_{\text{max}}$ and TCV do not converge to the correct physics ($\kappa \approx 0$), however they still converge, but 
to values corresponding to large $\kappa$ material's physics. %km: finde ich gut
In Table~\ref{Sneddon2dlayered}, the numbers of GMRES iterations are moderate for $\nu=0.2$. For higher Poisson ratios, we observe high linear iteration numbers. The incompressibility and the mesh adaptivity seem to significantly impact the linear solver.
We observe the same effects for $\kappa=10^{-8}$ in Table~\ref{Sneddon2dlayeredAdaptive}.

  \begin{table}[htbp!]
   \scriptsize % \tiny
\centering
\caption{Sneddon's pressure-driven cavity layered. Average number of GMRES iterations ($\varnothing$lin) per Newton step ($\#$AS). Computations with $Q_2^c Q_1^c Q_1^c$ elements for different problem size, $\epsilon=h$ for three Poisson ratios. Quantities of interest: COD$_{\text{max}}$ and TCV and $\kappa =10^{-8}$. Geometrically refined mesh as depicted on the right in Figure~\ref{SneddonGeoLayered}.
}\label{Sneddon2dlayeredAdaptive}
\begin{tabular}{|l|l|l|r|r|r|r|l|r|} \hline
\multicolumn{8}{|c|}{Sneddon layered adaptive: robustness in $h$, $\lambda$, $\epsilon$; $\epsilon=h$, $\kappa=10^{-8}$; mixed}\\ \hline\hline
 \rowcolor{gray!50} \multicolumn{1}{|c|}{FE$(u,p,\varphi)$} & \multicolumn{1}{|c|}{$\nu$} & \multicolumn{1}{|c|}{$h$} & \multicolumn{1}{|c|}{$\#$DoFs} & \multicolumn{1}{|c|}{$\varnothing$lin} & \multicolumn{1}{|c|}{$\#$AS} & \multicolumn{1}{|c|}{COD$_{\max}$} & \multicolumn{1}{|c|}{TCV}\\
 \hline
  \hline 
    \rowcolor{yellow!5}$Q_2^c Q_1^c Q_1^c$ & 0.2 & $0.353$ & $257\,924$  & 10 & 3 & 0.00242526 & 0.0107193\\ \hline %4+1
  \rowcolor{yellow!10}$Q_2^c Q_1^c Q_1^c$ & 0.2 & $0.176$  & $263\,604$  & 20 & 3 & 0.00221789 & 0.0080340\\ \hline %4+2
  \rowcolor{yellow!15}$Q_2^c Q_1^c Q_1^c$ & 0.2 & $0.088$  & $282\,484$  & 18 & 3 & 0.00208683 & 0.0069646\\ \hline %4+3
  \rowcolor{yellow!20}$Q_2^c Q_1^c Q_1^c$ & 0.2 & $0.044$ & $350\,804$  & 29 & 6 & 0.00200814 & 0.0064862\\ \hline %4+4
    \rowcolor{yellow!25}$Q_2^c Q_1^c Q_1^c$ & 0.2 & $0.022$ & $610\,164$  & 26 & 4 & 0.00196329 & 0.0062530\\ \hline %4+5
        \rowcolor{yellow!30}$Q_2^c Q_1^c Q_1^c$ & 0.2 & $0.011$ & $1\,620\,244$  & 32 & 3 & 0.00193890 & 0.0061344 \\ \hline %4+6
            \rowcolor{yellow!35}$Q_2^c Q_1^c Q_1^c$ & 0.2 & $0.0055$ & $5\,606\,324$  & 40 & 3 & 0.00192609 & 0.0060733\\ \hline %4+7
 ref.~\cite{SneddLow69} & 0.2 &  &  &  & & 0.00192000 & 0.0060318\\ \hline 
 \rowcolor{blue!5}$Q_2^c Q_1^c Q_1^c$ & 0.4999 & $0.353$ & $257\,924$  & 40 &  2 & 0.00223914 & 0.0116630\\ \hline %4+1
  \rowcolor{blue!10}$Q_2^c Q_1^c Q_1^c$ & 0.4999 & $0.176$  & $263\,604$  & 74 & 5 & 0.00187365 & 0.0077192\\ \hline %4+2
  \rowcolor{blue!15}$Q_2^c Q_1^c Q_1^c$ & 0.4999 & $0.088$  & $282\,484$  & 229 & 4 & 0.00168693 & 0.0060788\\ \hline %4+3
  \rowcolor{blue!20}$Q_2^c Q_1^c Q_1^c$ & 0.4999 & $0.044$ & $350\,804$  & 511 & 4 & 0.00159278 & 0.0053537\\ \hline %4+4
    \rowcolor{blue!25}$Q_2^c Q_1^c Q_1^c$ & 0.4999 & $0.022$ & $610\,164$  & 601 & 6 & 0.00154436 & 0.0050158\\ \hline %4+5
        \rowcolor{blue!30}$Q_2^c Q_1^c Q_1^c$ & 0.4999 & $0.011$ & $1\,620\,244$  & 565 & 5 & 0.00151941 & 0.0048527 \\ \hline %4+6
        \rowcolor{blue!35}$Q_2^c Q_1^c Q_1^c$ & 0.4999 & $0.0055$ & $5\,606\,324$  & 641 & 5 & 0.00150668 & 0.00477248\\ \hline %4+7 still running
        ref.~\cite{SneddLow69} & 0.4999 &  &  &  &  & 0.00150019 & 0.0047130\\ \hline 
        \hline
        %%%%%%%%%%%%%%%%%%%%%%%%%%%%%%%%%%%%%%%%%%%%%  
  \rowcolor{Maroon!5}$Q_2^c Q_1^c Q_1^c$ & 0.5 & $0.353$ & $257\,924$  & 40 & 2 & 0.00223891 & 0.0116629\\ \hline %4+1
  \rowcolor{Maroon!10}$Q_2^c Q_1^c Q_1^c$ & 0.5 & $0.176$  & $263\,604$  & 73 & 5 & 0.00187338 & 0.0077187\\ \hline %4+2
  \rowcolor{Maroon!15}$Q_2^c Q_1^c Q_1^c$ & 0.5 & $0.088$  & $282\,484$  & 227 & 4 & 0.00168667 & 0.0060782\\ \hline %4+3
       ref.~\cite{SneddLow69} & 0.5 &  &  &  &  &  0.0015000 & 0.0047124\\ \hline 
\end{tabular} 
\end{table}

In Table~\ref{Sneddon2dlayeredAdaptive}, the numerical results of the same tests are given as in Table~\ref{Sneddon2dlayered} for $\kappa=10^{-8}$. Analogously to Table~\ref{Sneddon2dKappSmallCgforS}, Table~\ref{Sneddon2dlayeredAdaptiveCgforS} contains the numerical results for the Sneddon test layered for high Poisson ratios and small $\kappa$. In contrast to Table~\ref{Sneddon2dlayeredAdaptive}, we approximate $\hat{S}^{-1}$ with a CG solver which is preconditioned with AMG. 

  \begin{table}[htbp!]
  \scriptsize %  \tiny
\centering
\caption{Sneddon's pressure-driven cavity layered with $Q_2^c Q_1^c Q_1^c$ elements and $\epsilon=h$. Average number of GMRES iterations ($\varnothing$lin) per Newton step ($\#$AS). CG plus AMG is used for $(g(\tilde{\varphi})A_u)^{-1}$ and $\hat{S}^{-1}$: the average number of CG iterations for $(g(\tilde{\varphi})A_u)^{-1}$ is 38 for $\nu=0.4999$ and 36 for $\nu=0.5$. The average number of CG iterations for $\hat{S}^{-1}$ is 8 for $\nu=0.4999$ and 7 for $\nu=0.5$. Computations for different problem size, $\epsilon=h$ for three Poisson ratios. Quantities of interest: COD$_{\text{max}}$ and TCV and $\kappa =10^{-8}$. Geometrically refined mesh as depicted on the right in Figure~\ref{SneddonGeoLayered}.
}\label{Sneddon2dlayeredAdaptiveCgforS}
%\textbf{Sneddon layered adaptive: robustness in $h$, $\lambda$, $\epsilon$; $\epsilon=2h$, $\kappa=10^{-8}$; %mixed}\\
%\vspace{0.2cm}
%\renewcommand*{\arraystretch}{1.2}
\begin{tabular}{|l|r|r|r|r|r|l|r|} \hline
\multicolumn{8}{|c|}{Sneddon layered adaptive: robustness in $h$, $\lambda$, $\epsilon$; $\epsilon=h$, $\kappa=10^{-8}$; mixed; CG for two blocks}\\ \hline\hline
  \rowcolor{gray!50} \multicolumn{1}{|c|}{FE$(u,p,\varphi)$} & \multicolumn{1}{|c|}{$\nu$} & \multicolumn{1}{|c|}{$h$} & \multicolumn{1}{|c|}{$\#$DoFs} & \multicolumn{1}{|c|}{$\varnothing$lin} &  \multicolumn{1}{|c|}{$\#$AS} & \multicolumn{1}{|c|}{COD$_{\max}$} & \multicolumn{1}{|c|}{TCV}\\
 \hline\hline
 \rowcolor{blue!5}$Q_2^c Q_1^c Q_1^c$ & 0.4999 & $0.353$ & $257\,924$  & 40 & 2 & 0.00223914 & 0.0116630\\ \hline %4+1
  \rowcolor{blue!10}$Q_2^c Q_1^c Q_1^c$ & 0.4999 & $0.176$  & $263\,604$  & 70 & 4 & 0.00187365 & 0.0077192\\ \hline %4+2
  \rowcolor{blue!15}$Q_2^c Q_1^c Q_1^c$ & 0.4999 & $0.088$  & $282\,484$  & 165 & 4 & 0.00168693 & 0.0060788\\ \hline %4+3
  \rowcolor{blue!20}$Q_2^c Q_1^c Q_1^c$ & 0.4999 & $0.044$ & $350\,804$  & 153 & 4 & 0.00159278 & 0.0053537\\ \hline %4+4
    \rowcolor{blue!25}$Q_2^c Q_1^c Q_1^c$ & 0.4999 & $0.022$ & $610\,164$   & 145 & 5 & 0.00154436 & 0.0050158\\ \hline %4+5
        \rowcolor{blue!30}$Q_2^c Q_1^c Q_1^c$ & 0.4999 & $0.011$ & $1\,620\,244$   & 139 & 5 & 0.00151941 & 0.0048527\\ \hline %4+6
        \rowcolor{blue!35}$Q_2^c Q_1^c Q_1^c$ & 0.4999 & $0.0055$ & $5\,606\,324$  & 148 & 5 & 0.00150668 & 0.0047724\\ \hline %4+7 
        ref.~\cite{SneddLow69} & 0.4999 &  &  &  & & 0.00150019 & 0.0047130\\ \hline 
        \hline
  \rowcolor{red!5}$Q_2^c Q_1^c Q_1^c$ & 0.5 & $0.353$ & $257\,924$  & 40 &  2 & 0.00223891 & 0.0116629\\ \hline %4+1
  \rowcolor{red!10}$Q_2^c Q_1^c Q_1^c$ & 0.5 & $0.176$  & $263\,604$  & 70 & 4 & 0.00187338 & 0.0077187\\ \hline %4+2
  \rowcolor{red!15}$Q_2^c Q_1^c Q_1^c$ & 0.5 & $0.088$  & $282\,484$  & 200 &  4 & 0.00168667 & 0.0060782\\ \hline %4+3
    \rowcolor{red!20}$Q_2^c Q_1^c Q_1^c$ & 0.5 & $0.044$ & $350\,804$  & 229 & 4 & 0.00159254 & 0.0060782\\ \hline %4+4
    \rowcolor{red!25}$Q_2^c Q_1^c Q_1^c$ & 0.5 & $0.022$ & $610\,164$   & 238  & 7 & 0.00154414 & 0.0050151\\ \hline %4+5
        \rowcolor{red!30}$Q_2^c Q_1^c Q_1^c$ & 0.5 & $0.011$ & $1\,620\,244$   & 220 & 5 & 0.00151920 & 0.0048521\\ \hline %4+6
        \rowcolor{red!35}$Q_2^c Q_1^c Q_1^c$ & 0.5 & $0.0055$ & $5\,606\,324$  & 222 & 5 & 0.00150648 & 0.0047718\\ \hline %4+7 
   ref.~\cite{SneddLow69} & 0.5 & &  &   &  &  0.00150000 & 0.0047124\\ \hline  
\end{tabular} 
\end{table}

The results of COD$_{\text{max}}$ and TCV in Tables~\ref{Sneddon2dlayeredAdaptive} and~\ref{Sneddon2dlayeredAdaptiveCgforS} look promising for all three Poisson ratios. For $\nu=0.5$ the solver does not converge with sufficiently small $\kappa$ and $h\to 0$. An explanation is given in Remark~\ref{kappa_problems} (Section~\ref{sec_sneddon}). In Table~\ref{Sneddon2dlayeredAdaptiveCgforS} for high Poisson ratios, the modified approximation of $\hat{S}^{-1}$ changes the behavior of the linear solver. With a relative tolerance of $10^{-6}$ for the preconditioned CG solver for $(g(\tilde{\varphi})A_u)^{-1}$ and $\hat{S}^{-1}$, we observe that more GMRES iterations are required. The number of linear iterations is relatively high, but nearly constant for $\nu=0.4999$ and $\nu=0.5$. The number of linear iterations increases for higher Poisson ratios with adaptive refined meshes and $\epsilon=h$. The results of COD$_{\text{max}}$ and TCV match the manufactured reference values.
% from Table~\ref{reference_cod_tcv}.

\begin{figure}[htbp!]
    \centering
    \includegraphics[width=0.3\textwidth]{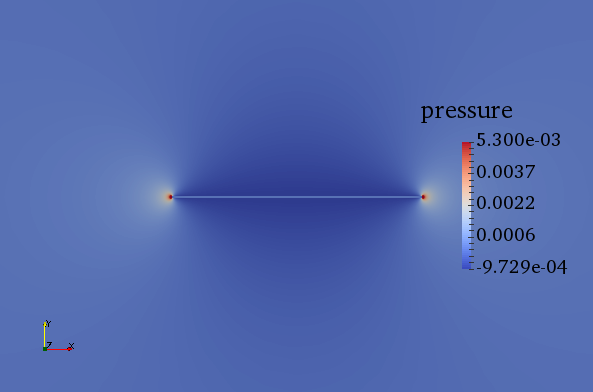}
    \hspace{0.25cm}
    \includegraphics[width=0.3\textwidth]{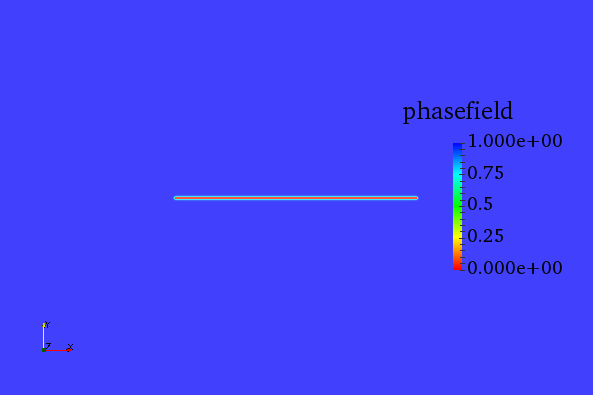}\\
    \vspace{1em}
    \includegraphics[width=0.3\textwidth]{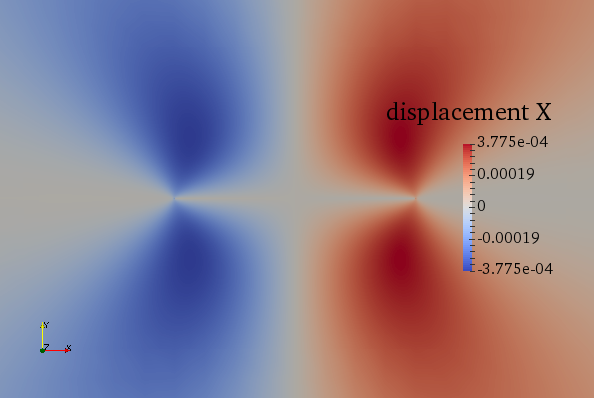}
    \hspace{0.25cm}
    \includegraphics[width=0.3\textwidth]{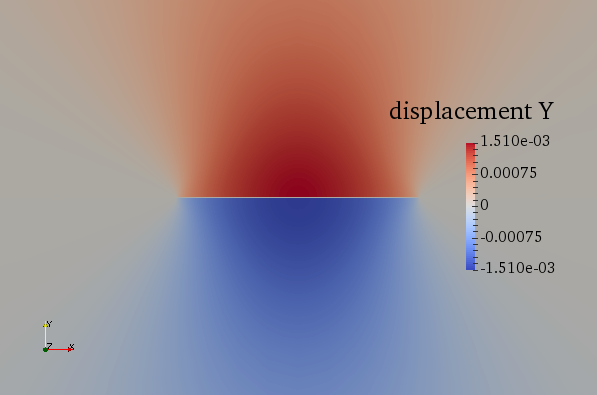}\\
    \vspace{1em}
    \includegraphics[width=0.3\textwidth]{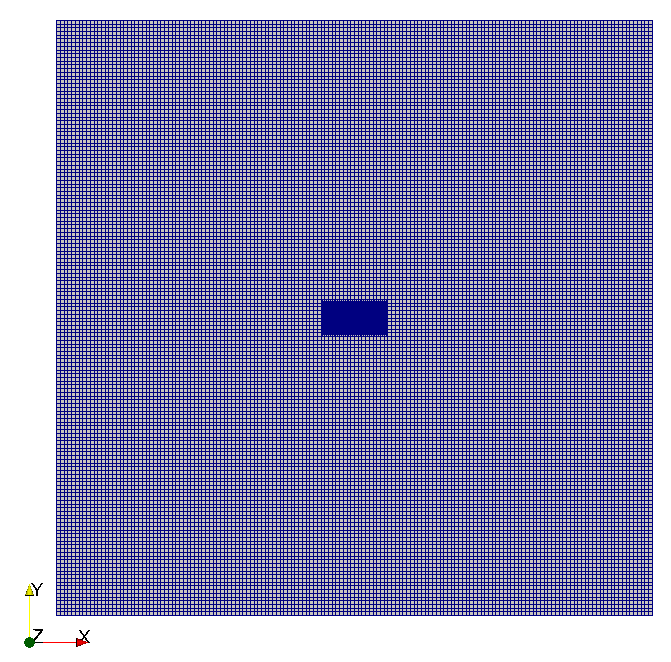}
    \hspace{0.25cm}
    \includegraphics[width=0.3\textwidth]{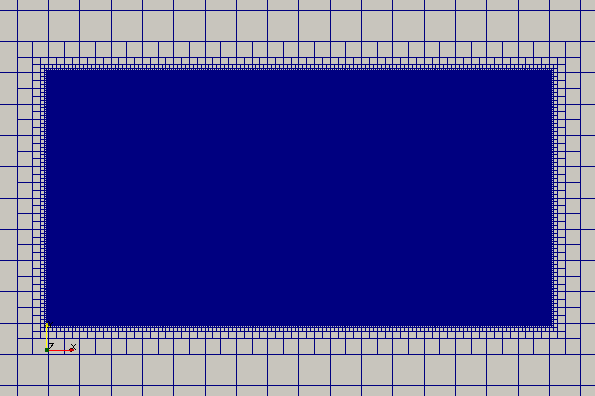}
    \caption{Sneddon 2d layered. Upper four snapshots: zoom-in solutions from left to right, and top to bottom: the pressure field, the phase-field, the displacement in the $x$-direction and the displacement in the $y$-direction. The solutions are for $\nu=0.5$ from Table~\ref{Sneddon2dlayeredAdaptiveCgforS} on the finest level with $Q_2^c Q_1^c Q_1^c$ elements. The solutions fit the reference values from~\cite{basava2020adaptive}.
    Lower two snapshots show on the left the whole domain $(-20,20)^2$ with the adaptively refined mesh in the last refinement step. On the bottom right, a zoom-in snapshot of the crack zone is given on the domain, where the upper snapshots are taken.}
    \label{snapshots_layered2dCGforS}
\end{figure}

In Figure~\ref{snapshots_layered2dCGforS}, the solutions of $u_x$, $u_y$, $p$, and $\varphi$ are presented as zoom-in snapshots for $\nu=0.5$ with a compressible layer, based on Table~\ref{Sneddon2dlayeredAdaptiveCgforS}. Especially the pressure field (upper left snapshot) is expected to have zero values in the interior of the crack and the maximal values in the crack tip on the left and the right of the pre-defined initial crack. Further, in Figure~\ref{snapshots_layered2dCGforS}, the mesh on the finest refinement level is given on the bottom left. On the bottom right, the crack zone is shown, on which the computed solutions are presented above to get an impression of the mesh size around the fracture.

%%%%%%%%%%%%%%%%%%%%%%%%%%%%%%%%%%%%%%%%%%%%%%%%%%%%%%%%%%%%%%%%%%%%%%%%%%%%%%%%%%%%%%%%%

\subsection{Single edge notched pure tension test}\label{sec_sent}

As the last example, we use the single-edge notched tension test from Miehe et al.~\cite{miehe2010phase} testing with three Poisson ratios. We use the predictor-corrector scheme from Heister et al.~\cite{heister2015primal} for two steps of adaptive mesh refinement on four times uniformly refined mesh with a phase-field threshold of $0.5$. The parameter setting is the same as in~\cite{miehe2010phase} but we use the mixed problem formulation and discretization from
Section~\ref{notation_equations} and vary the Poisson ratio; see Table~\ref{tension_nu_table}.
\begin{table}[htbp!]
  \scriptsize %  \tiny
\centering
\caption{Parameter setting for three tests with different Poisson's ratios for the single-edge notched tension test with $\kappa=10^{-8}$, and $\epsilon=4\,h$. The maximal number of DoFs is given in the last column for the test cases. For all tests, four uniform ($h=0.011$) and two adaptive refinement steps are conducted with a phase-field threshold of $0.5$ for predictor-corrector.}\label{tension_nu_table}
\begin{tabular}{|l|c|r|c|}\hline
 \rowcolor{gray!50}\multicolumn{1}{|c|}{$\nu$} & \multicolumn{1}{|c|}{$\mu$} & \multicolumn{1}{|c|}{$\lambda$} & \multicolumn{1}{|c|}{$\#$DoFs} \\ \hline \hline
$0.3$  & $80.77\cdot 10^3$ & $121.15\cdot 10^3$ & $19\,584$\\ \hline
$0.45$ & $80.77\cdot 10^3$ & $726.93\cdot 10^3$ & $19\,704$\\ \hline
$0.49$ & $80.77\cdot 10^3$ & $3957.73\cdot 10^3$ & $19\,498$\\ \hline
 \end{tabular}
\end{table}

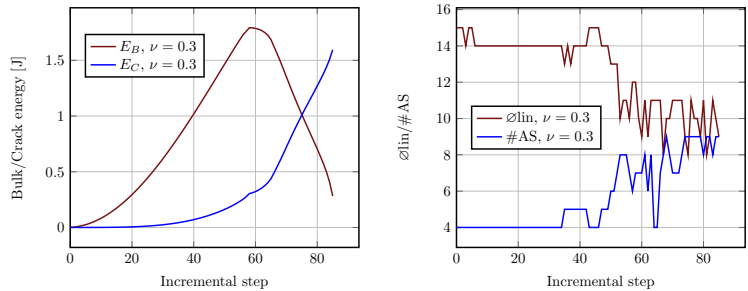
\begin{figure}[htbp!]
\centering
\begin{tikzpicture}[xscale=0.56,yscale=0.56]
\begin{axis}[
    ylabel = Bulk/Crack energy $\lbrack\mathrm{J}\rbrack$,
    xlabel = Incremental step,
     legend cell align=left,
    legend style={at={(0.05,0.8)},anchor=west},
 %legend pos=north west, 
 grid =major,
  xmin=0,
  ]
\addplot[Maroon]
table[x=step,y=bulk,col sep=comma] {Data/statistics_02_sent.csv}; 
\addlegendentry{$E_B$, $\nu=0.3$}
\addplot[blue] 
table[x=step,y=crack,col sep=comma] {Data/statistics_02_sent.csv};
\addlegendentry{$E_C$, $\nu=0.3$}
\end{axis}
\end{tikzpicture}
\hspace{0.2cm}
\begin{tikzpicture}[xscale=0.56,yscale=0.56]
\begin{axis}[
    ylabel = $\varnothing$lin/$\#$AS,
    xlabel = Incremental step,
     legend cell align=left,
 %legend pos=north west, 
     %ytick={0,20,40,60,80},
 legend style={at={(0.05,0.5)},anchor=west},
 grid =major,
  xmin=0,
  ]
\addplot[Maroon]
table[x=step,y=lin,col sep=comma] {Data/statistics_02_sent.csv}; 
\addlegendentry{$\varnothing$lin, $\nu=0.3$}
\addplot[blue] 
table[x=step,y=as,col sep=comma] {Data/statistics_02_sent.csv};
\addlegendentry{$\#$AS, $\nu=0.3$}
\end{axis}
\end{tikzpicture}
\caption{Left: bulk ($E_B$) and crack energy ($E_C$) for the single-edge notched tension test, AT$_2$ functional, adaptive refined meshes. The incremental step size was $10^{-4}\,\mathrm{s}$ for the first 58 steps and reduced to $10^{-5}\,\mathrm{s}$. Right: number of linear iterations on average per Newton step ($\varnothing$lin), and number of active set/Newton steps ($\#$AS) against the incremental steps. The crack starts propagating at incremental step 66, mesh refinement starts at step 57. $\nu=0.3$.}\label{bulk_crack_sent_03}
\end{figure}

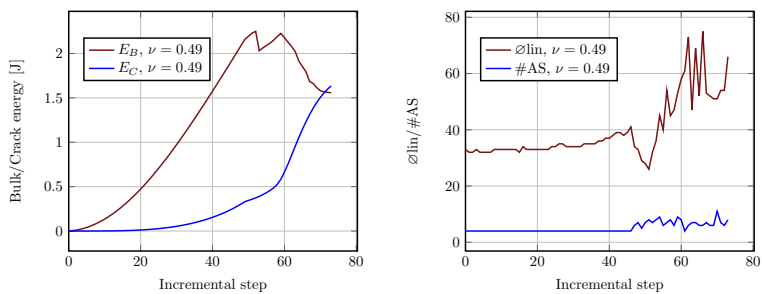
\begin{figure}[htbp!]
\centering
\begin{tikzpicture}[xscale=0.56,yscale=0.56]
\begin{axis}[
    ylabel = Bulk/Crack energy $\lbrack\mathrm{J}\rbrack$,
    xlabel = Incremental step,
 %legend pos=north west,
  legend cell align=left,
 legend style={at={(0.05,0.8)},anchor=west},
 grid =major,
  xmin=0,
  ]
\addplot[Maroon]
table[x=step,y=bulk,col sep=comma] {Data/statistics_049_sent.csv}; 
\addlegendentry{$E_B$, $\nu=0.49$}
\addplot[blue] 
table[x=step,y=crack,col sep=comma] {Data/statistics_049_sent.csv};
\addlegendentry{$E_C$, $\nu=0.49$}
\end{axis}
\end{tikzpicture}
\hspace{0.2cm}
\begin{tikzpicture}[xscale=0.56,yscale=0.56]
\begin{axis}[
 legend cell align=left,
    ylabel = $\varnothing$lin/$\#$AS,
    xlabel = Incremental step,
 legend pos=north west, 
 legend style={at={(0.05,0.8)},anchor=west},
 grid =major,
  xmin=0,
  ]
\addplot[Maroon]
table[x=step,y=lin,col sep=comma] {Data/statistics_049_sent.csv}; 
\addlegendentry{$\varnothing$lin, $\nu=0.49$}
\addplot[blue] 
table[x=step,y=as,col sep=comma] {Data/statistics_049_sent.csv};
\addlegendentry{$\#$AS, $\nu=0.49$}
\end{axis}
\end{tikzpicture}
\caption{Left: bulk ($E_B$) and crack energy ($E_C$) for the single-edge notched tension test, AT$_2$ functional, adaptive refined meshes. The incremental step size was $10^{-4}\,\mathrm{s}$ for the first 48 steps and reduced to $10^{-5}\,\mathrm{s}$. Right: number of linear iterations on average per Newton step ($\varnothing$lin), and number of active set/Newton steps ($\#$AS) against the incremental steps. The crack starts propagating at incremental step 60, mesh refinement starts at step 52. $\nu=0.49$.}\label{bulk_crack_sent_049}
\end{figure}

We consider the bulk and crack energy as two further numerical quantities of interest.
The bulk energy $E_B$ can be computed via
\begin{align*}
E_B(u,\varphi)= \int_{\Omega} (g(\tilde{\varphi}) \psi(E_{\text{lin}}(u))\, \mathrm{d}{(x,y)},
\end{align*}
where the strain energy functional is defined as
\begin{align*}
 \psi(E_{\text{lin}}(u))\coloneqq \mu \tr \left(E_{\text{lin}}(u)^2\right) + \frac{1}{2} \lambda \tr \left(E_{\text{lin}}(u)\right)^2.
\end{align*}
Here, no manufactured reference values are provided and we only present
values computed numerically.
Further, we compute the crack energy $E_C$ via
\begin{align*}
 E_C(u,\varphi) = \frac{G_C}{2} \int_{\Omega} \left( \frac{(\varphi-1)^2}{\epsilon} +  \epsilon |\nabla \varphi|^2\right) \, \mathrm{d}{(x,y)}.
\end{align*}
Again, no manufactured reference values are provided. At least for $\nu=0.3$, we can compare our results for $E_B$ and $E_C$ with reference values from the literature, e.g.,~\cite{AmGeLo15,mang2020phase}. In Figures~\ref{bulk_crack_sent_03} and~\ref{bulk_crack_sent_049}, on the left side, the bulk and the crack energy are plotted versus the incremental step number. On the right of Figures~\ref{bulk_crack_sent_03} to~\ref{bulk_crack_sent_049}, the average number of linear iterations and the number of Newton/AS steps are plotted. 
The number of linear iterations behaves differently for $\nu=0.3$ from the results for higher Poisson ratios. While for $\nu=0.3$, in Figure~\ref{bulk_crack_sent_03} on the right, the linear iterations decrease if the crack starts propagating, in Figure~\ref{bulk_crack_sent_049}, the linear iterations increase up to an average of more than 70 iterations at the end of the crack simulations. 
%Mathematically it is not trivial if we have a solution on a Neumann boundary (as on the left boundary in this test) where the crack propagates to. 
%%%%%%%%%%%%%%%%%%%%%%%%%%%%%%%%%%%%%%%%%%%%%%%%%%

\begin{figure}[htbp!]
    \centering
    \includegraphics[height=3.05cm]{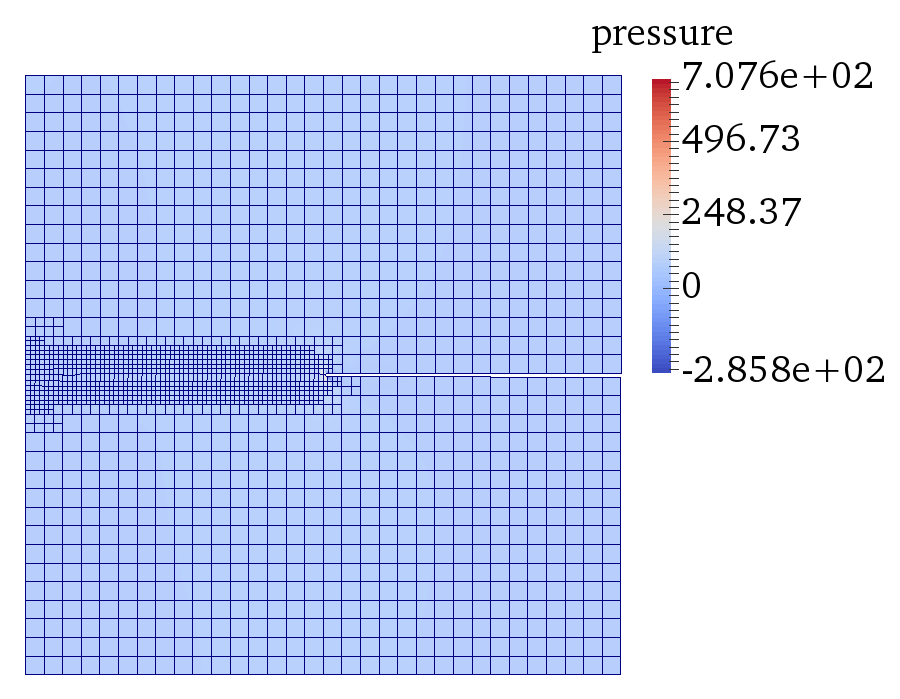}
    \hspace{0.37cm}
    \includegraphics[height=3.05cm]{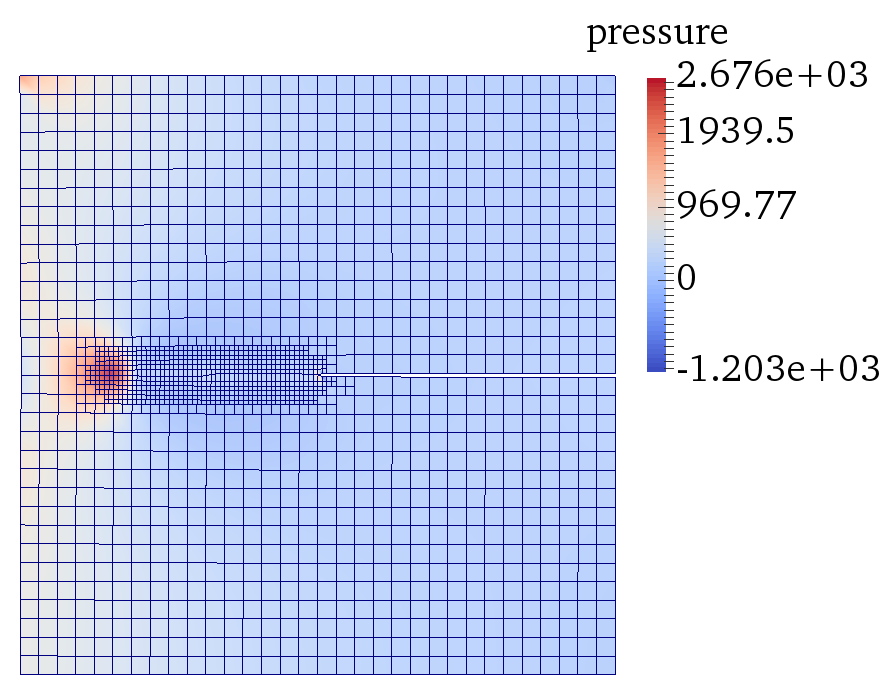}
        \hspace{0.28cm}
    \includegraphics[height=3.05cm]{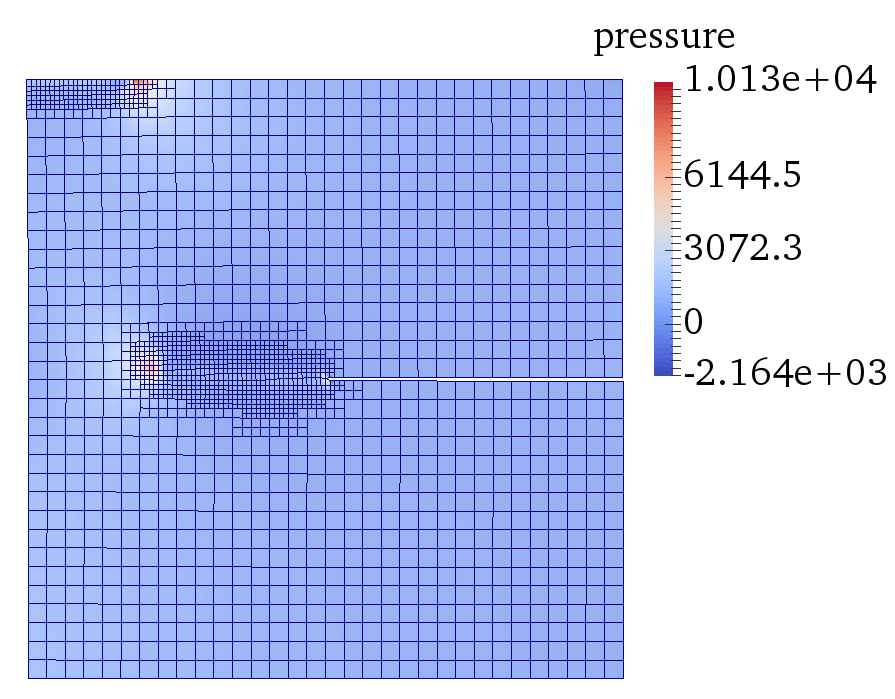}\\
    \includegraphics[height=3.1cm]{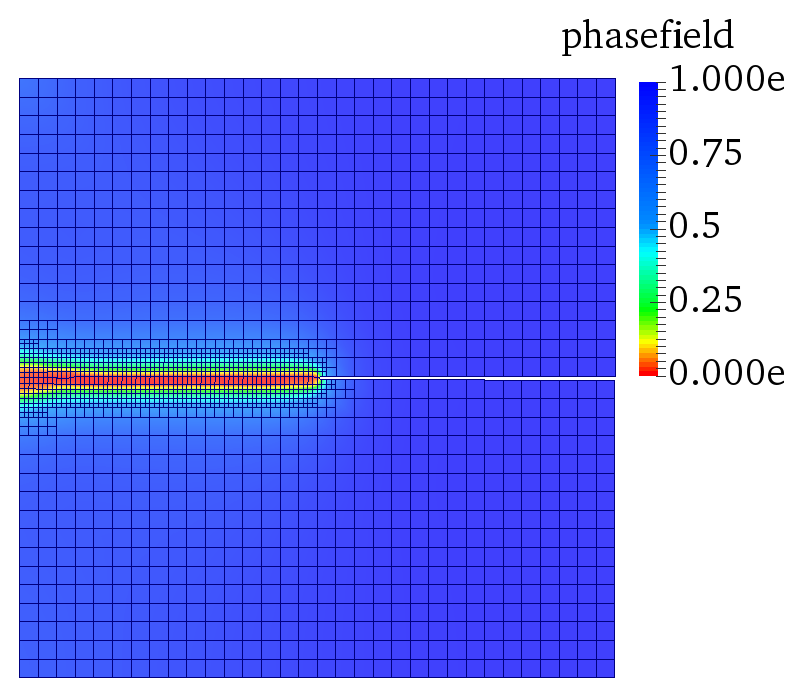}
        \hspace{0.7cm}
    \includegraphics[height=3.1cm]{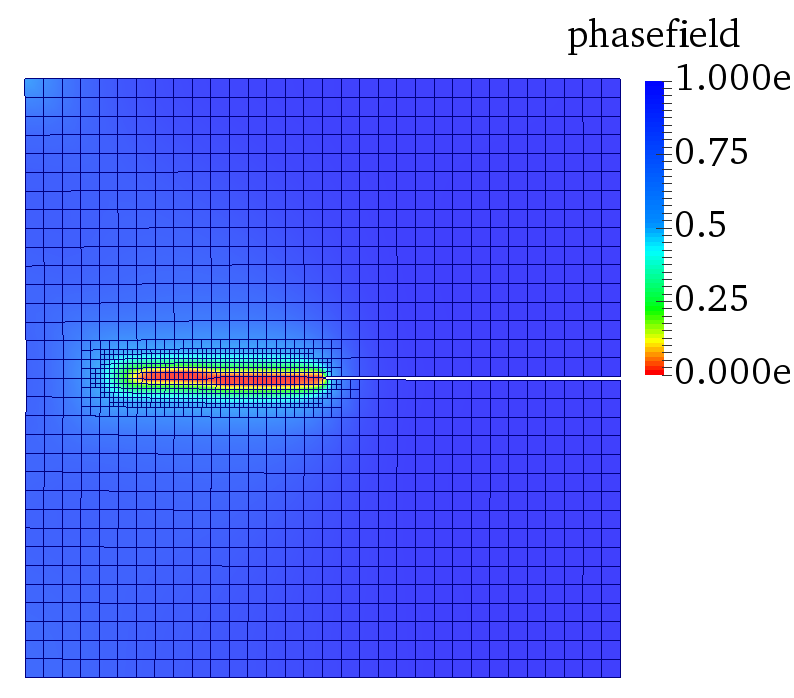}
        \hspace{0.7cm}
    \includegraphics[height=3.1cm]{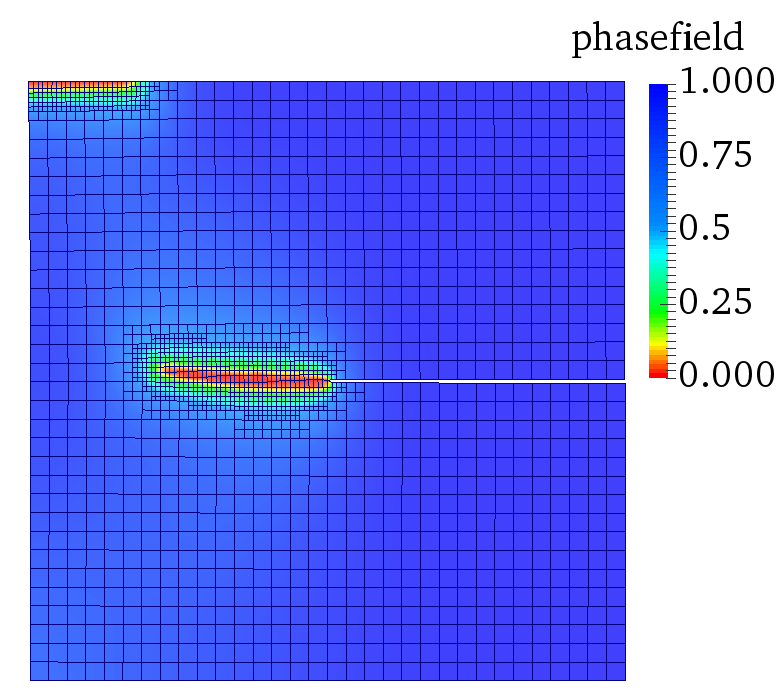}
    \caption{Snapshots of the solution for the single-edge notched tension test with $\nu=0.3, 0.45,$ and $0.49$ from left to right. Adaptive mesh refinement with predictor-corrector. Pressure field (first row) and phase-field (second row) on deformed two-dimensional domain in incremental step 88 for $\nu=0.3$, step 70 for $\nu=0.45$, and step 75 for $\nu=0.49$.}
    \label{snapshots_sent_049}
\end{figure}

In Figure~\ref{snapshots_sent_049}, snapshots of the pressure field and phase-field are given for $\nu=0.49$, where -- to the author's knowledge -- no reference values are available in the literature. The crack paths look similar as for $\nu=0.3$, but a slight asymmetry is visible in the crack path. We decided to present the crack path during the simulation to depict the pressure field with the maximal value in front of the crack tip while the pressure values in the crack are zero.
%%%%%%%%%%%%%%%%%%%%%%%%%%%%%%%%%%%%%%%%%%%%%%%%%%%%%%%%%%%%%%%%
The computed bulk and crack energies in Figure~\ref{bulk_crack_sent_03} fit well to results in the literature, e.g.,~\cite{heister2015primal}. The bulk energy increases until the critical energy release rate is reached, and the crack energy increases when the crack propagates while the bulk energy releases. Also, in Figure~\ref{bulk_crack_sent_049}, the bulk and crack energy curves fit the observed crack pattern in Figure~\ref{snapshots_sent_049}. For $\nu=0.49$ with snapshots in the last column in Figure~\ref{snapshots_sent_049}, no comparable results in the literature are available. The crack pattern differs from the snapshots for smaller Poisson ratios. We observe that the crack has an orientation to the upper left corner, and a second crack develops from the singularity in the corner, where non-homogeneous Dirichlet boundary conditions and Neumann boundary conditions meet. 
In the first column in Figure~\ref{snapshots_sent_049}, the pressure and phase-field solution is given for $\nu=0.3$ after total failure. The crack propagates from the center of the geometry to the left boundary, as we expect it. Further, one can see a pure zero pressure field after total failure. 

%%%%%%%%%%%%%%%%%%%%%%%%%%%%%%%%%%%%%%%%%%%%%%%%%%%%%%%%%%%%%%%%
\section{Conclusions}
\label{sec_conclusions}
In this work, a preconditioner
for a mixed formulation phase-field fracture model
that is
robust in $h,\epsilon$, and $\lambda$
 was developed and tested on four 
numerical examples for different Poisson ratios up to the incompressible limit, namely 
$\nu\to 0.5$ yielding $\lambda\to\infty$.
For the first test case, a hanging block with a slit, we confirmed the robustness and efficiency of the physics-based preconditioner, discretized with $Q_2^c Q_1^c Q_1^c$ finite elements. To the best of the authors' knowledge, in the last test case the well-known single edge notched tension test was considered for higher Poisson ratios for the first time. For $\nu=0.49$ a non-symmetric crack behavior and crack initiation from the upper left corner singularity was observed.
In Sneddon's test case and $\kappa=10^{-8}$, an impact of $\kappa$ on the condition of the $\kappa$-dependent block entries of the system matrix could be explicitly seen.

It is well-known that from a phase-field perspective the regularization parameter $\epsilon$ 
is challenging, in particular its choice in relation to $h$. However, we found in this 
paper that from a preconditioner perspective the first regularization parameter $\kappa$ (in the 
bulk term of the displacement equation) causes difficulties instead. 
Basically, we deal with 
an elliptic (Laplacian) term 
where diffusion ranges from $\kappa \approx 10^{-8}$ in the crack region to 1, a
difference of 8 orders of magnitude.
We expect that a carefully designed geometric multigrid preconditioner or a weighted BFBT preconditioner might  handle this situation better. 
%Both are major extensions out of scope in this work.
We emphasize that having robustness in $h,\epsilon$ and $\lambda$ is a significant 
contribution, which has not yet been studied so far in the published literature.
A second future extension would be thermodynamically consistent constitutive materials laws,
namely incorporating stress splitting in $\sigma(u,p)$.

\section*{Acknowledgements}
   K. Mang thanks Clemson University for the financial support for a one-month research stay.
   
    T. Heister was partially supported by the National Science Foundation (NSF)
Award DMS-2028346, OAC-2015848, EAR-1925575, by the Computational
Infrastructure in Geodynamics initiative (CIG), through the NSF under Award
EAR-0949446 and EAR-1550901 and The University of California -- Davis, and by
Technical Data Analysis, Inc. through US Navy STTR Contract N68335-18-C-0011.

Clemson University is acknowledged for generous allotment of compute time on Palmetto cluster.

      The second and third authors were supported by the German Research Foundation, Priority Program 1748 (DFG SPP 1748) within the subproject 
      Structure Preserving Adaptive Enriched Galerkin Methods for Pressure
      3D Fracture Phase-Field Models (WI 4367/2-1) with the project number 
      392587580.

%%%%%%%%%%%%%%%%%%%%%%%%%%%%%%%%%%%%%%%%%%%%%%%%%%
%\section*{References}
\bibliographystyle{plain}
\bibliography{references}
%%%%%%%%%%%%%%%%%%%%%%%%%%%%%%%%%%%%%%%%%%%%%%%%%%
\end{document}